\documentclass[12pt,a4paper]{article}

\usepackage{amsfonts}
\usepackage{mathrsfs}
\usepackage{indentfirst, latexsym, bm}
\usepackage{amsfonts,amssymb,amsthm}
\usepackage{amsbsy}
\usepackage{latexsym}
\usepackage[mathscr]{eucal}
\usepackage{amsmath,amscd}  

\usepackage{graphicx}

\newtheorem{theorem}{\hspace{1.3em}Theorem}[section]
\newtheorem{lemma}{\hspace{1.3em}Lemma}[section]

\newtheorem{corollary}{\hspace{1.3em}Corollary}[section]

\newtheorem{example}{\hspace{1.3em}Example}[section]

\begin{document}

\title{$k$-flaw Preference Sets}

\author{Po-Yi Huang$^{a,}$\thanks{Partially supported by NSC 96-2115-M-006-012
}
 \and Jun Ma$^{b,}$\thanks{Email address of the corresponding author: majun@math.sinica.edu.tw}
 \and  Jean Yeh$^{c,}$\thanks{jean.yh@ms45.url.com.tw}}

\date{}
\maketitle \vspace{-1cm} \begin{center} \footnotesize
 $^{a}$ Department of Mathematics, National Cheng Kung University, Tainan, Taiwan\\
$^{b}$ Institute of Mathematics, Academia Sinica, Taipei, Taiwan\\
$^{c}$ Department of Mathematics, National Taiwan University,
Taipei, Taiwan
\end{center}
\thispagestyle{empty}\vspace*{.4cm}

\begin{abstract}
In this paper, let $\mathcal{P}_{n;\leq s;k}^l$ denote a set of
$k$-flaw preference sets $(a_1,\ldots,a_n)$ with $n$ parking spaces
satisfying that $1\leq a_i\leq s$ for any $i$ and $a_1=l$ and
$p_{n;\leq s;k}^l=|\mathcal{P}_{n;\leq s;k}^l|$. We use a
combinatorial approach to the enumeration of $k$-flaw preference
sets by their leading terms. The approach relies on bijections
between the $k$-flaw preference sets and labeled rooted forests.
Some bijective results between certain sets of $k$-flaw preference
sets of distinct leading terms are also given. We derive some
formulas and recurrence relations for the sequences $p_{n;\leq
s;k}^l$ and give the generating functions for these sequences.
\end{abstract}


\newpage

\section{Introduction}
Throughout the paper, let $[n]:=\{1,2,\ldots,n\}$ and
$[m,n]:=\{m,\ldots, n\}$. Suppose that $n$ cars have to be parked in
$m$ parking spaces which are arranged in a line and numbered $1$ to
$m$ from left to right. Each car has initial parking preference
$a_i$; if space $a_i$ is occupied, the car moves to the first
unoccupied space to the right. We call $(a_1,\ldots,a_n)$ {\it
preference set}. Clearly, the number of preference sets is $m^n$. If
a preference set $(a_1,\ldots,a_n)$ satisfies $a_i\leq a_{i+1}$ for
$1\leq i\leq n-1$, then we say that the preference set is {\it
ordered}. If all the cars can find a parking space, then we say the
preference set is a {\it parking function}. If there are exactly $k$
cars which can't be parked, then the preference set is called a {\it
$k$-flaw preference set}.

Let $n$, $m$, $s$, and $k$ be four nonnegative integers with $1\leq
s\leq m$ and $k\leq n-1$. Suppose that  there are $m$ parking
spaces. We use $\mathcal{P}_{n,m;\leq s;k}$ to denote a set of
$k$-flaw preference sets $(a_1,\ldots,a_n)$ of length $n$ satisfying
$1\leq a_i\leq s$ for all $i$. For $1\leq l\leq s$, we use
$\mathcal{P}_{n,m;\leq s;k}^l$ to denote a set of preference sets
$(a_1,\ldots ,a_n)\in\mathcal{P}_{n,m;\leq s;k}$ such that $a_1=l$.
Let $p_{n,m;\leq s;k}=|\mathcal{P}_{n,m;\leq s;k}|$ and $p_{n,m;\leq
s;k}^l=|\mathcal{P}_{n,m;\leq s;k}^l|$. For any of the above cases,
if the parameter $k$ ( resp. $m$ ) doesn't appear, we understand
$k=0$ ( resp. $m=n$
 ); if the parameter $m$ and $s$ are both erased, we understand
$s=m=n$.

There are some results about parking functions with $s=m=n$. Riordan
introduced parking functions in \cite{R}. He derived that the number
of parking functions of length $n$ is $(n+1)^{n-1}$, which coincides
with the number of labeled trees on $n+1$ vertices by Cayley's
formula. Several bijections between the two sets are known (e.g.,
see \cite{FR,R,SMP}). Furthermore, define a generating function
$P(x)=\sum\limits_{n\geq 0}\frac{(n+1)^{n-1}}{n!}x^n$. It is well
known that $xP(x)$ is the compositional inverse of the function
$\psi(x)=xe^{-x}$, i.e., $\psi(xP(x))=x$. Riordan concluded that the
number of ordered parking functions is
$\frac{1}{n+1}{2n\choose{n}}$, which is also equals the number of
Dyck path of semilength $n$. Parking functions have been found in
connection to many other combinatorial structures such as acyclic
mappings, polytopes, non-crossing partitions, non-nesting
partitions, hyperplane arrangements,etc. Refer to
\cite{F,FR,GK,PS,SRP,SRP2} for more information.

Any parking function $(a_1,\ldots,a_n)$ can be redefined such that
its increasing rearrangement $(b_1,\ldots,b_n)$ satisfies $b_i\leq
i$.
 Pitman and  Stanley generalized the notion of parking functions
in \cite{PS}. Let ${\bf x}=(x_1,\ldots,x_n)$ be a sequence of
positive integers. The sequence $\alpha=(a_1,\ldots,a_n)$  is called
an ${\bf x}$-parking function if the non-decreasing rearrangement
$(b_1,\ldots,b_n)$ of $\alpha$ satisfies $b_i\leq x_1+\ldots +x_i$
for any $1\leq i\leq n$. Thus, the ordinary parking function is the
case ${\bf x}=(1,\ldots,1)$. By the determinant formula of
Gon\v{c}arove polynomials, Kung and Yan \cite{KY} obtained the
number of ${\bf x}$-parking functions for an arbitrary ${\bf x}$.
See also \cite{Y1,Y2,Y3} for the explicit formulas and properties
for some specified cases of ${\bf x}$.

An ${\bf x}$-parking function $(a_1,\ldots,a_n)$ is said to be
$k$-leading if $a_1=k$. Let $q_{n,k}$ denote the number of
$k$-leading ordinary parking functions of length $n$. Foata and
Riordan \cite{FR} derived a generating function for $q_{n,k}$
algebraically. Recently, Sen-peng Eu, Tung-shan Fu and Chun-Ju Lai
\cite{EFL} gave a combinatorial approach to the enumeration of
$(a,b,\ldots,b)$-parking functions by their leading terms.

Riordan \cite{R} told us the relations between ordered parking
functions and Dyck paths. Sen-peng Eu et al. \cite{EFY,ELY}
considered the problem of the enumerations of lattice paths with
flaws. It is natural to consider the problem of the enumerations of
preference sets with flaws. There is a interesting facts. Salmon and
Cayley in 1849 established the classical configuration of $27$ lines
in a general cubic surface. Given a line $l$, the number of lines,
which are disjoint from, intersect or are equal to  $l$, are
$16$,$10$ and $1$, respectively, see \cite{henderson1911} for the
detail information. These number exactly are the number of $0$-,
$1$- and $2$-flaw preference sets of length $3$. Peter J Cameron et
al. \cite{Cam} indicate that there are some relations between
$k$-flaw preference sets and the context of hashing since data would
be lost. Also they counted the number of $k$-flaw preference sets
and calculate the asymptotic. Ordered $k$-flaw preference sets were
studied in \cite{H1} and some enumerations for some parking
functions were given in \cite{H2}.

In this paper, we use the methods developed by Sen-peng Eu et al.
\cite{EFL} to study $k$-flaw preference sets. Sen-peng Eu et al.
find that {\it triplet-labelled rooted forests} enable $(a, 1,
\ldots, 1)$-parking functions to be manipulated on forests easily.
The methods are different with that in \cite{Cam}. We consider
triplet-labelled rooted forests associated with $k$-flaw preference
sets and establish a bijection between $k$-flaw preference sets and
label rooted forests, so as to enumerate $k$-flaw preference sets by
leading term.

First, we enumerate $k$-flaw preference sets in the set
$\mathcal{P}_{n;\leq s;k}$. Then we consider the triplet-labelled
rooted forests and bijections associated with the parking functions
in the set $\mathcal{P}_{n,n+k;\leq n+k}$. Using these bijections,
we find that $p_{n,n+k;\leq n+k}^1=p_{n,n+k;\leq n+k}^2=\ldots
=p_{n,n+k;\leq n+k}^{k+1}$ for any $k\geq 0$. Taking $l=1$, we have
(1) $p^1_{n,n+k;\leq n+k}=p_{n-1,n+k;\leq n+k}$ for any $k\geq 0$
and $n\geq 1$; (2) $p_{n,n+k, \leq n+k}^1
=\sum\limits_{i=1}^{n}{n-1\choose{i-1}}p_{i,i+k-1,\leq
i+k-1}^1p_{n-i}$ for any $k\geq 1$. When $k+1\leq l\leq n+k-1$, we
obtain the recurrence relation $p_{n,n+k;\leq n+k}^{l}-p_{n,n+k;\leq
n+k}^{l+1}={n-1\choose{l-k-1}}p_{l-k-1,l-1;\leq l-1}p_{n+k-l}.$

To enumerate $k$-flaw preference sets in the set
$\mathcal{P}_{n;\leq s;k}^l$, we study the triplet-labelled rooted
forests and bijections associated with the $k$-flaw preference sets
in the set $\mathcal{P}_{n;\leq s;k}$. We prove that $p_{n;\leq
s;k}^1=p_{n;\leq s;k}^2=\ldots =p_{n;\leq s;k}^k$ for any $k\geq 1$;
letting $l=1$, we have $p_{n;\leq
s;k}^1=\sum\limits_{i=1}^{s-k-1}{n-1\choose{s-i-k-1}}p_{n+k-s+i;\leq
i}p_{s-k-i,s-i-1;\leq s-i-1}^1$ for any $k\geq 1$ and $k+1\leq s\leq
n$. For the cases with $s=l$, we derive some interesting identities:
(1) $p_{n;\leq s;k+1}^s=p_{n;\leq s;k}^1$ for any $k\geq 1$ and
$k+2\leq s\leq n$; (2) $p_{n+1;\leq s;k}^s=p_{n;\leq s;k}$ for
$k\geq 0$ and $k+1\leq s\leq n$; (3) $p_{n+1;\leq s;k}^s=p_{n;\leq
s;k}$ for any $k\geq 0$ and $k+1\leq s\leq n$; (4)
$p_{n;1}^n=p_{n}^2$ for any $n\geq 2$; (5)
$p_{n;1}^{n}-p_{n;1}^{n-1}=p_{n-2}$ for any $n\geq 2$. When $k\leq
l\leq s-1$, we obtain the following recurrence relation
\begin{eqnarray*}p_{n;\leq s;k}^{l+1}-p_{n;\leq
s;k}^{l}&=&\sum\limits_{i=k-1}^{l-2}{n-1\choose{i-k+1}}p_{i-k+1,i;\leq
i}\left[ p_{n+k-i-1;\leq s-i-1}^{l-i}-p_{n+k-i-1;\leq
s-i-1}^{l-i-1}\right]\\
&&+{n-1\choose{l-k}}p_{l-k,l-1;\leq l-1}\left[ p_{n+k-l;\leq
s-l}^{1}-p_{n+k-l;\leq s-l;1}^{s-l}\right].\end{eqnarray*}

We also are interested in generating functions for some sequences.
For any $k\geq 1$ and $s\geq 0$, define a generating function
$D_{k,s}(x)=\sum\limits_{n\geq s+k}\frac{p_{n;\leq n-s;k}}{n!}x^n$,
then $D_{k,s}(x)=
[P(x)]^{k+1}\sum\limits_{i=0}^{k+s}\frac{(-1)^i(k+s+1-i)^i}{i!}x^i-[P(x)]^{k}\sum\limits_{i=0}^{k+s-1}\frac{(-1)^i(k+s-i)^i}{i!}x^i.$
Furthermore, let $D(x,y,z)=\sum\limits_{k\geq 1}\sum\limits_{s\geq
0}\sum\limits_{n\geq s+k}\frac{p_{n;\leq n-s;k}}{n!}x^ny^sz^k$, then
$D(x,y,z)
=\frac{zP(x)}{y-zP(x)}\left[\frac{P(x)-y}{e^{xy}-y}-\frac{(1-z)P(x)}{e^{xzP(x)}-zP(x)}\right].
$

For any $l\geq 0$ and $k\geq 0$, define a generating function
$H_{l,k}(x)=\sum\limits_{n\geq l+1}\frac{p_{n,n+k;\leq
n+k}^{n+k-l}}{(n-1)!}x^n$, then $H_{l,k}(x)$ satisfies the
recurrence relation $H_{l,k}(x)=H_{l-1,k}(x)-\frac{p_{l,l+k;\leq
l+k}^{1}}{(l-1)!}x^l+\frac{p_{l}}{l!}x^{l+1}[P(x)]^{k+1}$ with the
initial conditions $H_{0,k}(x) =x[P(x)]^{k+1}. $ Let
$H_{k}(x,y)=\sum\limits_{l\geq 0}H_{l,k}(x)y^l$ and
$H(x,y,z)=\sum\limits_{k\geq 0}H_{k}(x,y)z^k$ , then
$H_k(x,y)=\frac{xP(xy)\{[P(x)]^{k+1}-y[P(xy)]^{k+1}\}}{1-y}$ and
$H(x,y,z)=\frac{xP(xy)}{1-y}\left[\frac{P(x)}{1-zP(x)}-\frac{yP(xy)}{1-zP(xy)}\right].$

Define a generating function $W(x,y,z,v)=\sum\limits_{k\geq
1}\sum\limits_{s\geq 0}\sum\limits_{l\geq s}\sum\limits_{n\geq
k+l}\frac{p_{n;\leq n-s;k}^{n-l}}{(n-1)!}x^ny^lz^sv^k,$ then
 \begin{eqnarray*}
W(x,y,z,v)&=&\frac{xyvP(xy)}{y-1}\left\{\left[\frac{yP(xy)}{yz-vP(xy)}-\frac{P(x)}{yz-vP(x)}\right]R(xy,z)\right.\\
&&+\left.\frac{P(x)}{yz-vP(x)}R(xy,\frac{v}{y}P(x))-\frac{yP(xy)}{yz-vP(xy)}R(xy,\frac{v}{y}P(xy))\right\}\\
&&+\frac{vP(x)}{yz-vP(x)}\left[F(x,y,z)-F(x,y,\frac{v}{y}P(x))\right]\end{eqnarray*}
where $R(x,y)=\frac{P(x)-y}{e^{xy}-y}$ and
$F(x,y,z)=\frac{x}{e^{xyz}-z}\left[\frac{P(xy)(P(x)-yP(xy))}{1-y}-\frac{zP(xyz)[P(x)-yzP(xyz)]}{1-yz}\right].$

Recently, Postnikov and Shapiro \cite{postnikov2004} gave a new
generalization, building on work of Cori, Rossin and Salvy
\cite{cori2002}, the $G$-parking functions of a graph. For the
complete graph $G=K_{n+1}$, the defined functions in
\cite{postnikov2004} are exactly the classical parking functions.
So, in the future work, we will consider $k$-flaw $G$-parking
function.

We organize this paper as follows. In Section $2$, we enumerate
$k$-flaw preference sets in the set $\mathcal{P}_{n;\leq s;k}$. In
Section $3$, we consider the triplet-labelled rooted forests and
bijections associated with the parking functions in the set
$\mathcal{P}_{n,n+k;\leq n+k}$. In Section $4$, we give the
enumerations of parking functions in the set
$\mathcal{P}_{n,n+k;\leq n+k}^l$. In Section $5$, we study the
triplet-labelled rooted forests and bijections associated with the
$k$-flaw preference sets in the set $\mathcal{P}_{n;\leq s;k}$. In
Section 6, we investigate the problems of the enumerations of
preference sets in the set $\mathcal{P}_{n;\leq s;k}^l$. In Section
$7$, we obtain some generating functions for some sequences given in
the previous sections. In Appendix, we list the values of $p_{n;\leq
s;k}^l$ for $n\leq 7$ and $p_{n,n+k;\leq n+k}^l$ for any $n\leq 5$
and $k\leq 3$.

\section{Counting the number of elements in $\mathcal{P}_{n;\leq s;k}$}
In this section, we will consider the enumerations of preference
sets in the set $\mathcal{P}_{n;\leq s;k}$.
\begin{lemma}\label{lemmap(n;<=s;k)com}
Let $1\leq s\leq n$ and $1\leq k\leq s-1$, then
$$p_{n;\leq
s;k}=\sum\limits_{i=1}^{s-k}{n\choose{n-s+i+k}}p_{s-i-k,s-i-1;\leq
s-i-1}p_{n-s+i+k;\leq i}$$
\end{lemma}
{\bf Proof.} For any $\alpha=(a_1,\ldots,a_n)\in\mathcal{P}_{n;\leq
s;k}$,
 suppose that the $(n-i)$-th parking
space is the last empty one. Obviously, $n-s+1\leq i\leq n-k$. Let
$S=\{j\mid a_j>n-i\}$ and $\alpha_S$ be a subsequence of $\alpha$
determined by the subscripts in $S$. Then $|S|=k+i$. Let
$T=[n]\setminus S$, then $|T|=n-k-i$ and $a_j<n-i$ for any $j\in T$.
Let $\alpha_T$ be a subsequence of $\alpha$ determined by the
subscripts in $T$, then $\alpha_T\in\mathcal{P}_{n-k-i,n-i-1;\leq
n-i-1}$. Suppose $\alpha_s=(b_1,\ldots,b_{k+i})$, then
$(b_1-n+i,\ldots,b_{k+i}-n+i)\in\mathcal{P}_{i+k;\leq s+i-n}$.

There are ${n\choose{i+k}}$ ways to choose $i+k$ numbers from $[n]$
for the elements in $S$. There are $p_{i+k;\leq s+i-n}$ and
$p_{n-k-i,n-i-1;\leq n-i-1}$ possibilities for $\alpha_S$ and
$\alpha_T$, respectively. Hence, for any $1\leq k\leq s-1$, we have,
\begin{eqnarray*}p_{n;\leq s;k}&=&\sum\limits_{i=n-s+1}^{n-k}{n\choose{i+k}}p_{n-i-k,n-i-1;\leq
n-i-1}p_{i+k;\leq s+i-n}\\
&=&\sum\limits_{i=1}^{s-k}{n\choose{n-s+i+k}}p_{s-i-k,s-i-1;\leq
s-i-1}p_{n-s+i+k;\leq i}.\end{eqnarray*}\hfill$\blacksquare$

\section{Triple-labelled rooted forests and a bijection}
\addtocounter{figure}{1} In this section,  we consider the
triplet-labelled rooted forests and bijections associated with the
parking functions in the set $\mathcal{P}_{n,n+k;\leq n+k}$. Using
these bijections, we give the enumerations of parking functions in
the set $\mathcal{P}_{n,n+k;\leq n+k}^l$.

 Let $\mathcal{B}_{n,k}$
be a set of all sequences $(T_0,\ldots,T_{k})$ of length $k+1$ such
that (1) the union of the vertex sets of $T_0,\ldots,T_{k}$ is
$\{R_i\mid 0\leq i\leq k\}\cup [n]$, where $R_i\notin [n]$ is just
an artificial label; (2) each $T_i$ is a tree with root $R_i$; (3)
$T_i$ and $T_j$ are disjoint if $i\neq j$.

Let $F\in \mathcal{B}_{n,k}$. For any $x\in [n]$, there is an unique
root $R_i$ which is connected with $x$. Define the height of $x$ to
be the number of edges connecting $x$ with root $R_i$. If the height
of a vertex $z$ is less than the height of $x$ and $\{z,x\}$ is an
edge of $F$, then $z$ is the predecessor of $x$, $x$ is a child of
$z$, and write $z={\rm pre}(x)$ and $x\in {\rm child}(z)$.

Fixing a sequence $F$ of rooted trees in $\mathcal{B}_{n,k}$, we
define a linear order $<_F$ on $[n]$ by the following rules. Let
$x,y\in [n]$.

(1) For any $i\neq j$, $x\in T_i$ and $y\in T_j$, if $i<j$, then
$x<_Fy$.

(2) For any $i$ and $x,y\in T_i$, if the height of $x$ is less than
the height of $y$, then $x<_Fy$; if the height of $x$ is equals the
height of $y$, and ${\rm pre}(x)<_F{\rm pre}(y)$, then $x<_Fy$.

(3) For any $i$ and $x,y\in T_i$, if ${\rm pre}(x)={\rm pre}(y)$,
but $x<y$, then $x<_Fy$.

The sequence formed by writing $\{1,\ldots,n\}$ in the increasing
order with respect to $<_F$ is denoted by
$\sigma_F^{-1}=(\sigma^{-1}_F(1),\ldots,\sigma^{-1}_F(n))$. And the
permutation $\sigma_F$ is the inversion of $\sigma_F^{-1}$.

Next, we define the {\it forest specification} of $F$. Let $m_i$ be
the number of the vertices in $\bigcup\limits_{j=0}^{i-1}T_j$ for
$1\leq i\leq k+1$. Clearly, $m_{k+1}=n+k+1$. Set ${\bf
r}_F=(r_1,\ldots,r_{n+k})$ as follows.

(1) $r_1$ is the number of children of the vertex $R_0$.

(2) $r_i$ is the number of children of the vertex
$\sigma^{-1}_F(i-1)$ if $2\leq i\leq m_1$.

(3) $r_{m_i+1}$ is the number of children of the vertex $R_{i}$ for
$1\leq i\leq k$.

(4) $r_i$ is the number of children of the vertex
$\sigma^{-1}_F(i-j)$ if $m_{j-1}+2\leq i\leq m_{j}$ for some $j$.

Given $F\in \mathcal{B}_{n,k}$, we may obtain $({\bf
r}_F,\sigma_F)$. Let
$$a_i=\left\{\begin{array}{lll}
1&\text{if}&1\leq\sigma_F(i)\leq r_1\\
k+1&\text{if}&1+\sum\limits_{i=1}^{k}r_i\leq\sigma_F(i)\leq
\sum\limits_{i=1}^{k+1}r_i\end{array}\right.$$ It is well known that
each labeled tree on $n+1$ vertices would correspond to a parking
function of length $n$. So, we obtain a parking function
$\alpha_F=(a_1,\ldots,a_n)\in\mathcal{P}_{n,n+k;\leq n+k}$ and
parking spaces $m_1,\ldots,m_{k}$ could not be occupied.

Conversely, for any $\alpha=(a_1,\ldots,a_n)\in\mathcal
{P}_{n,n+k;n+k}$, let ${\bf r}_\alpha=(r_1,\ldots,r_{n+k})$ be the
$specificatioin$ of $\alpha$, i.e., $r_i=|\{j\mid a_j=i\}|.$
Furthermore, we may suppose that parking spaces $m_1,\ldots,m_k$ are
empty, then $r_{m_i}=0$ for $i=1,\ldots,k$.

For $1\leq i\leq n$, define
$$\pi_{\alpha}(i)=|\{a_j\mid \text{ either }a_j<a_i,\text{ or }
a_j=a_i \text{ and }j<i\}|.$$ Note that
$\pi_{\alpha}=(\pi_{\alpha}(1),\ldots,\pi_{\alpha}(n))$ is a
permutation of $[n]$. Let $\pi_{\alpha}^{-1}$ be the inversion of
$\pi_{\alpha}$.

Let $$t_i=\left\{\begin{array}{lll}
m_1-1&\text{if}&i=0\\
m_{i+1}-m_{i}-1&\text{if}&i=1,\ldots,k-1\\
n+k-m_k&\text{if}&i=k\\
\end{array}\right.$$

By $m_1,\ldots,m_k$, we can decompose $\bf r_\alpha$ and
$\pi^{-1}_\alpha$ into the following $k+1$
subsequences,respectively:

$${\bf r_{i}}=(r_{i,1},\ldots,r_{i,t_i})=\left\{\begin{array}{lll}
(r_1,\ldots,r_{m_1-1})&\text{if}&i=0\\
(r_{m_{i}+1},\ldots,r_{m_{i+1}-1)}&\text{if}&i=1,\ldots,k-1\\
(r_{m_{k}+1},\ldots,r_{n+k)}&\text{if}&i=k+1\\
\end{array}\right.$$
and
$$\sigma_{i}=(\sigma_i(1),\ldots,\sigma_i(t_i))=\left\{\begin{array}{lll}
(\pi^{-1}(1),\ldots,\pi^{-1}(m_1-1))&\text{if}&i=0\\
(\pi^{-1}(m_{i}-i+1),\ldots,\pi^{-1}(m_{i+1}-i-1))&\text{if}&i=1,\ldots,k-1\\
(\pi^{-1}(m_k-k+1),\ldots,\pi^{-1}(n))&\text{if}&i=k+1\\
\end{array}\right.$$

 We associate $({\bf r_i},\sigma_i)$ with  a rooted
tree $T_{i}$ on $t_i+1$ vertices. The vertex set of $T_{i}$ is
$\{R_{i}\}\cup\{\sigma_{i}(j)\mid 1\leq j\leq t_{i}\}$, where
$R_i\notin [n]$ is just an artificial label. Let $R_{i}$ be the root
of $T_{i}$. The children of $R_{i}$ are
$\sigma_{i}(1),\ldots,\sigma_{i}(r_{i,1})$. For any $1\leq j\leq
t_{i}$, the children of $\sigma_{i}(j)$ are
$\sigma_{i}(1+\sum\limits_{m=0}^jr_{i,m}),\ldots,\sigma_{i}(\sum\limits_{m=0}^{j+1}r_{i,m})$.
So, we obtain a rooted forest $F_\alpha$ with $k+1$ components
$T_0,\ldots,T_k$ and the vertex set of $F_\alpha$ is $\{R_i\mid
0\leq i\leq k\}\cup[n]$.

\begin{lemma}\label{bijectionphi}
There is a bijection $\phi$ between $\mathcal{P}_{n,n+k;\leq n+k}$
and $\mathcal{B}_{n,k}$.
\end{lemma}

\begin{corollary}\label{tripleroot}
Let $\alpha=(a_1,\ldots,a_n)\in\mathcal{P}_{n,n+k;\leq n+k}$ and
$F=(T_0,\ldots,T_k)=\phi(\alpha)\in \mathcal{B}_{n,k}$. Suppose
$x\in [n]$, $x$ is a child of root $R_s$ for some $0\leq s\leq k$.
Let $\mu$ be the number of the non-root vertices of
$\bigcup\limits_{j=0}^{s-1}T_j$, then $a_x=\mu+s+1$.
\end{corollary}
{\bf Proof.} Since $\alpha\in\mathcal{P}_{n,n+k;\leq n+k}$, we could
obtain $({\bf r}_\alpha,\pi_\alpha^{-1})$ and suppose that parking
spaces $m_1$,$\ldots$,$m_k$ couldn't be occupied. Observe that
$a_j=m_i+1$ for some $j\in [n]$ if and only if $j$ is a child of
$R_i$ in $F$. Hence, we have $a_x=m_s+1$. So $a_x=\mu+s+1$ since
$\mu=m_s-s$.\hfill$\blacksquare$

\begin{corollary}\label{triplenonroot}
Let $\alpha\in\mathcal{P}_{n,n+k;\leq n+k}$ and
$F=(T_0,\ldots,T_k)=\phi(\alpha)\in \mathcal{B}_{n,k}$. Suppose
$x_1,x_2\in [n]$, $x_1x_2$ is an edge of $T_s$ for some $0\leq s\leq
k$ and $\pi_\alpha(x_1)<\pi_\alpha(x_2)$, then
$a_{x_2}=\pi_{\alpha}(x_1)+s+1$.
\end{corollary}

{\bf Proof.} Since $\alpha\in\mathcal{P}_{n,n+k;\leq n+k}$, we could
obtain $({\bf r}_\alpha,\pi_\alpha^{-1})$ and suppose that parking
spaces $m_1$,$\ldots$,$m_k$ couldn't be occupied. Then there are the
subsequences $({\bf r_s},\sigma_s)$ of $({\bf
r}_\alpha,\pi_\alpha^{-1})$ such that $x_1=\sigma_s(j)$ and
$x_2=\sigma_s(l+\sum\limits_{m=0}^jr_{i,m})$ for some $j$ and $l$.
Note that $\pi_\alpha(x_1)=m_{s}-s+j$ and $a_{x_2}=m_{s}+j+1$.
Hence,
$a_{x_2}=\pi_\alpha(x_1)+s+1$.\hfill$\blacksquare$\\

By Corollary \ref{tripleroot} and Corollary \ref{triplenonroot}, we
may associate $\alpha$ with a ($k+1$)-component rooted forest
$\hat{F}_{\alpha}$ on $n+k+1$ vertices, called
$triplet$-$labeled~rooted ~forest$. Since
$\alpha\in\mathcal{P}_{n,n+k;\leq n+k}$, we suppose that parking
spaces $m_1,\ldots,m_k$ couldn't be occupied. Let $m_0=0$ and
$\mu_i=m_i-i$ for any $0\leq i\leq k$. The vertex set of
$\hat{F}_{\alpha}$ is $\{(R_i,0,\mu_i)\mid 0\leq i\leq
k\}\cup\{(i,a_i,\pi_{\alpha}(i))\mid 1\leq i\leq n\}$ of $triplets$,
where $R_i\notin [n]$ is just an artificial label for discriminating
the additional triplets. Let $(R_0,0,\mu_0),\ldots,(R_k,0,\mu_{k})$
be the root of distinct trees of $\hat{F}_{\alpha}$. For any two
vertices $u=(x_1,y_1,z_1)$ and $v=(x_2,y_2,z_2)$, $v$ is a child of
$u$ if there exists $i$ such that $\mu_i\leq z_1,z_2\leq\mu_{i+1}$
and $y_2=z_1+i+1$.

For example, take $n=12$ and
$\alpha=(7,1,11,12,6,11,12,2,6,2,6,12,1)\in\mathcal{P}_{12,14;\leq
14}$. It is easy to check that $m_1=4$ and $m_2=11$. Hence,
$\mu_0=0$, $\mu_1=3$ and $\mu_2=9$. We can obtain the following
table:
$$\begin{array}{|l|r|r|r|r|r|r|r|r|r|r|r|r|}
\hline
i&1&2&3&4&5&6&7&8&9&10&11&12\\
a_i&6&1&8&12&7&12&5&8&12&2&1&5\\
\pi_\alpha(i)&6&1&8&10&7&11&4&9&12&3&2&5\\
\hline
\pi^{-1}_\alpha(i)&2&11&10&7&12&1&5&3&8&4&6&9\\
\hline
\end{array}
$$
and ${\bf r}_\alpha=(2,1,0,0,2,1,1,2,0,0,0,3,0,0)$. The rooted
forest $F_\alpha$ associated with $\alpha$ is shown on
Fig.\thefigure.
\begin{center}
\includegraphics[width=3in]{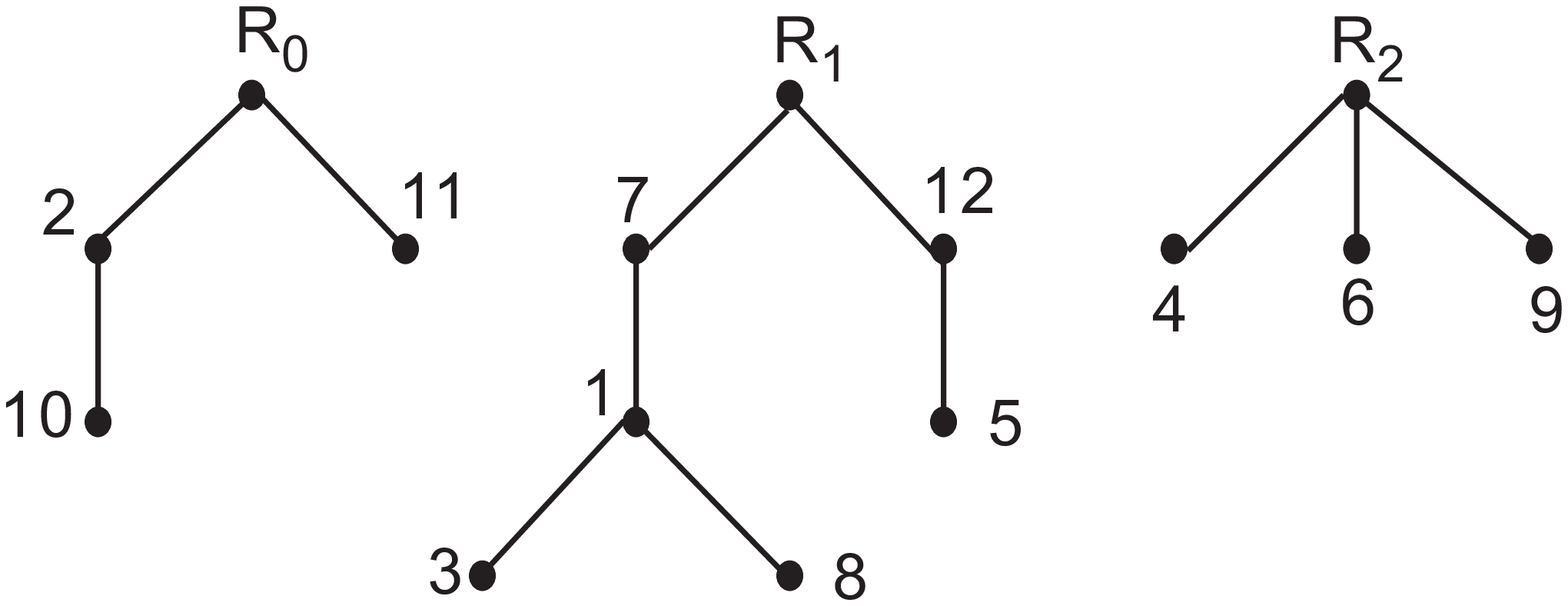}\\
Fig.\thefigure. A rooted forest $F_\alpha$ associated with the
parking function $\alpha=(7,1,11,12,6,11,12,2,6,2,6,12,1)$.
\end{center}
\addtocounter{figure}{1}

The triplet-labelled rooted forest $\hat{F}_\alpha$ associated with
$\alpha$ is shown on Fig.\thefigure.
\begin{center}
\includegraphics[width=3.75in]{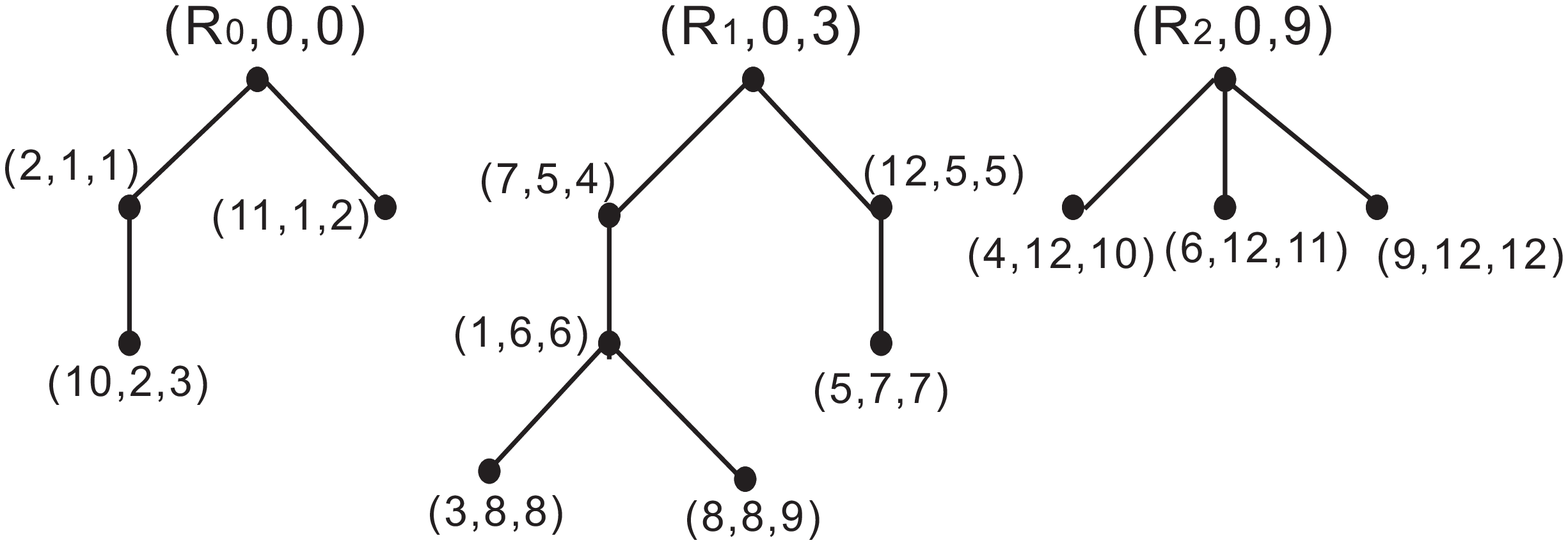}\\
Fig.\thefigure. A triplet-labelled rooted forest $\hat{F}_\alpha$
associated with the parking function
$\alpha=(7,1,11,12,6,11,12,2,6,2,6,12,1)$.
\end{center}

For any $\alpha\in\mathcal{P}_{n,n+k;\leq n+k}$, by triplet-labelled
rooted
 forest $\hat{F}_\alpha$ associated with $\alpha$, we could
 obtain the bijection $\phi$ between $\mathcal{P}_{n,n+k;\leq n+k}$ and
$\mathcal{B}_{n,k}$.

Given $\alpha\in\mathcal{P}_{n,n+k;\leq n+k}$, let $\hat{F}_\alpha$
be a triplet-labelled rooted forest $\hat{F}_\alpha$ associated with
$\alpha$. Let $\psi(x,y,z)=x$ for any $(x,y,z)\in\{(R_i,0,\mu_i)\mid
0\leq i\leq k\}\cup\{(i,a_i,\pi_{\alpha}(i))\mid 1\leq i\leq n\}$,
and two vertices $\psi(x_1,y_1,z_1)$ and $\psi(x_2,y_2,z_2)$ are
adjacent if and only if $(x_1,y_1,z_1)$ and $(x_2,y_2,z_2)$ are
adjacent in $\hat{F}$. So, we get a labeled forest $F_\alpha$. Let
$\phi(\alpha)=F_\alpha$, clearly, $\phi(\alpha)\in
\mathcal{B}_{n,k}$.

To describe $\phi^{-1}$, for each $F\in \mathcal{B}_{n,k}$, let
$F=(T_0,\ldots,T_k)$ and the tree $T_i$ have the root $R_i$ for each
$i$. We express $F$ in a form, called {\it canonical form}, of a
plane rooted forest. Let $T_0$, $\ldots$, $T_k$ be placed from left
to right. If a vertex has more than one child then the labels of
these children are increasing from left to right. Let $\mu_0=0$ and
$\mu_i$ be the number of the non-root vertices in
$\bigcup\limits_{j=0}^{i-1}T_j$ for $1\leq i\leq k$. Let
$\theta(R_i)=(R_i,0,\mu_i)$ for each $i$. For any $j\in [n]$, let
$\theta(j)=(j,y_j,z_j)$, where $y_j$ and $z_j$ are determined by the
following algorithm A.

{\bf Algorithm A.}

(1) Let $F=(T_0,\ldots,T_k)$. Traverse $T_i$ by a breadth-first
search, suppose $j$ is the $s$-th non-root vertex and label the
third entries $z_j=s+\mu_i$.

(2) For any two $v=(x_1,y_1,z_1)$ and $u=(x_2,y_2,z_2)$, if $u$ is a
child of $v$ and $x_1=R_i$, then $y_2=\mu_i+i+1$; if $u$ is a child
of $v$, the vertices $x_1$ and $x_2$ are in $T_i$ and $x_1,x_2\in
[n]$, then $y_2=z_1+i+1$.

Note that if $u$ is a child $v$ in $T_i$, then
$y_2-\mu_i-i=z_2-\mu_i=1$ if $x_1=R_i$; and $z_2>z_1$ implies that
$y_2-\mu_i-i=(z_1+i+1)-\mu_i-i\leq z_2-\mu_i$ if $x_1\neq R_i$.
Sorting the triplets of non-vertices by the first entries, the
sequence $\varphi^{-1}(F)=(y_1,\ldots,y_n)$, which is formed by
their second entries, is the required parking function. Obviously,
parking spaces $\mu_i+i$ couldn't be occupied for $1\leq
i \leq k$.\hfill$\blacksquare$\\

For any $\alpha=(a_1,\ldots,a_n)\in\mathcal{P}_{n,n+k;\leq n+k}$, we
may obtain ${\bf r}_\alpha=(r_1,\ldots,r_{n+k})$  and
$\pi_{\alpha}=(\pi_{\alpha}(1),\ldots,\pi_{\alpha}(n))$ and suppose
that parking spaces $m_1,\ldots,m_k$ couldn't be occupied with
$m_1<m_2<\ldots <m_k$, then define the following parameters about
$\alpha$:

(1) $l_\alpha$ : the leading term of $\alpha$, i.e., $l_\alpha
=a_1$;

(2) $n_\alpha$ : the number of $a_1$ in $\alpha$, i.e.,
$t_\alpha=|\{i\mid a_i=a_1\}|$;

(3) $m_{\alpha}=m_k$;

(4) $g_{\alpha}=max\{m_i\mid m_i<l_\alpha\}$;

(5) $h_{\alpha}=max\{i\mid m_i<l_\alpha\}$;

(6) $\tau_{\alpha}=|\{i\mid a_j<l_\alpha \}|+1$;

Clearly, $\tau_\alpha\geq l_{\alpha}-h_{\alpha}$. Now, let
$\mathcal{P}_{n,n+k;\leq n+k}^l=\{\alpha\in\mathcal{P}_{n,n+k;\leq
n+k}\mid l_\alpha=l\}$. Let $\mathcal{F}_{n,k}^l$ denote the set of
triplet-labelled rooted forests $F'_{\alpha}$ associated with
$\alpha\in\mathcal{P}_{n,n+k;\leq n+k}^l$. If $\mathcal{R}$ is a set
consisting of some parking functions, then we always use
$\mathcal{F}(\mathcal{R})$ to denote the set of triplet-labelled
rooted forests associated with the parking functions in
$\mathcal{R}$.

The following lemma has the same proof as the lemma in \cite{EFL}.
For the sake of completeness, we still prove it as follows.

\begin{lemma}\label{bijection<m1}
Let $k\geq 0$. For any $1\leq l\leq n+k-1$, let
$\mathcal{A}_1=\{\alpha\in\mathcal{P}_{n,n+k; \leq n+k}^l\mid
\tau_\alpha>l-h_\alpha \text{ or } n_\alpha\geq 2\}$ and
$\mathcal{C}_1=\{\alpha\in\mathcal{P}_{n,n+k; \leq n+k}^{l+1}\mid
l\geq g_\alpha+1\}$. Then there is a bijection from $\mathcal{A}_1$
to $\mathcal{C}_1$.
\end{lemma}
{\bf Proof.} It suffices to establish a bijection
$\psi:\mathcal{F}({\mathcal{A}_1})\rightarrow\mathcal{F}({\mathcal{C}_1})$.
Given an $\hat{F}_\alpha\in\mathcal{F}(\mathcal{A}_1)$, let
$u=(1,l,\pi_\alpha(1))\in \hat{F}_\alpha$. Obviously,  $u\in
T_{h_\alpha}$ and $l\geq g_\alpha+1$. If $\tau_\alpha>l-h_\alpha$,
then $\pi_\alpha(1)>l-h_\alpha$ and $T_{h_\alpha}-T_{h_\alpha}(u)$
has at least $l-h_\alpha$ vertices. On the other hand, $u\in
T_{h_\alpha}$ implies that there are at least $l-h_\alpha-1$ terms
$a_j$ satisfying $g_\alpha+1\leq a_j<l$. Hence, if $n_\alpha\geq 2$,
then there are at least $l-h_\alpha$ vertices in
$T_{h_\alpha}-T_{h_\alpha}(u)$ as well.

Let $(R_{h_\alpha},0,g_\alpha-{h_\alpha})$ be the $0$-th vertex.
Traverse $T_{h_\alpha}-T_{h_\alpha}(u)$ by breadth-first search and
locate the $(l-g_\alpha)$-th vertex in
$T_{h_\alpha}-T_{h_\alpha}(u)$, say $v$. By attaching
$T_{h_\alpha}(u)$ to $v$ so that $u$ is the first child of $v$ in
$T_{h_\alpha}-T_{h_\alpha}(u)$, updating the second and the third
entries of all non-root vertices by Algorithm A, and the other trees
remain unchangeable, we obtain $\psi(\hat{F}_\alpha)$. the triplet
of $u$ becomes $(1,l+1,\pi(1))$. $l\geq g_\alpha+1$ implies that
$\psi(\hat{F}_\alpha)\in\mathcal{F}(\mathcal{C}_1)$.

To find $\psi^{-1}$, given an
$\hat{F}_\beta\in\mathcal{F}(\mathcal{C}_1)$, let
$u=(1,l+1,\pi_\beta(1))\in \hat{F}_\beta$ and $v$ the parent of $u$.
Clearly, $u\in T_{h_\beta}$. Since $l\geq g_\beta+1$, $v$ isn't the
root $(R_{h_\beta},0,g_\beta-{h_\beta})$ of $T_{h_\beta}$. In
$\hat{F}_\beta$, we locate the vertex, say $w$, the third entry of
which is equal to $l-h_{\beta}-1$. Attach $T_{h_{\beta}}(u)$ to $w$
so that $u$ is the first child of $w$. By Algorithm A, the updated
triplet of $u$ becomes $(1,l,\pi(1))$ and the other trees are
unchangeable. We observe that either $\pi(1)=l$ if $v$ is another
child of $w$, or $\pi(1)>l$ otherwise. Hence
$\psi^{-1}(\hat{F}_\beta)\in\mathcal{F}(\mathcal{R})$.
\hfill$\blacksquare$\\
\addtocounter{figure}{1}

 For example, take $n=12$, $k=2$ and $l=2$.
We consider a parking function
$\alpha=(2,4,8,12,9,1,12,8,4,12,1,3)\in\mathcal{P}_{12,14;\leq 14}$.
Then $h_\alpha=0$, $n_\alpha=1$, $g_\alpha=0$ and $\tau_\alpha=3$.
Observe that $\tau_\alpha>l-h_\alpha$. On Fig.\thefigure~ is the
forest $\hat{F}_\alpha$ associated with $\alpha$.  Let $u=(1,2,3)$.
Note that $v=(11,1,2)$ is the $2$-th vertex of
$\hat{F}_\alpha-T_1(u)$ that is visited by a breadth-first search.
On Fig.$4$ is the corresponding forest $\psi(\hat{F}_\alpha)$, which
is obtained from $\hat{F}_\alpha-T_1(u)$ with $T_1(u)$ attached to
$v$ and with the second and third entries of the triplet updated.
Sorting the triplets of non-root vertices by the first entries, we
retrieve the corresponding parking function
$\alpha=(3,4,8,12,9,1,12,8,4,12,1,3)$ with leading term $3$ from
their second entries.
\begin{center}
\includegraphics[width=3.75in]{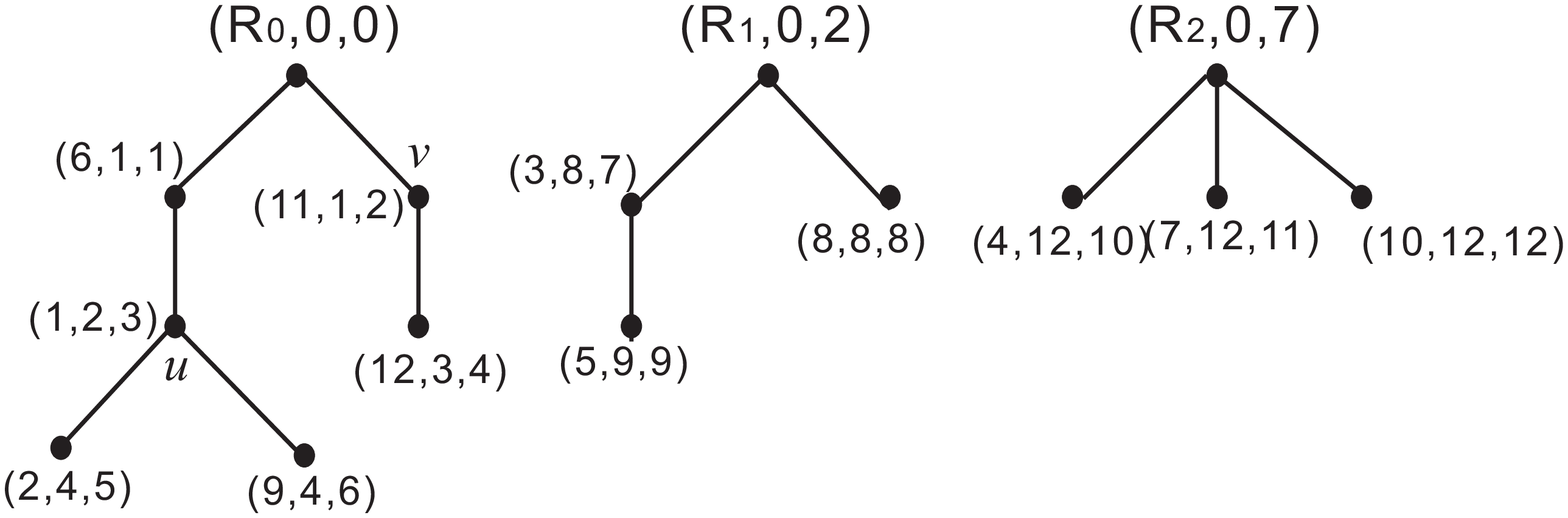}\\
Fig.\thefigure. The forest $F'_\alpha$associated with the parking
function $\alpha=(2,4,8,12,9,1,12,8,4,12,1,3)$
\end{center}
\addtocounter{figure}{1}
\begin{center}
\includegraphics[width=3.9in]{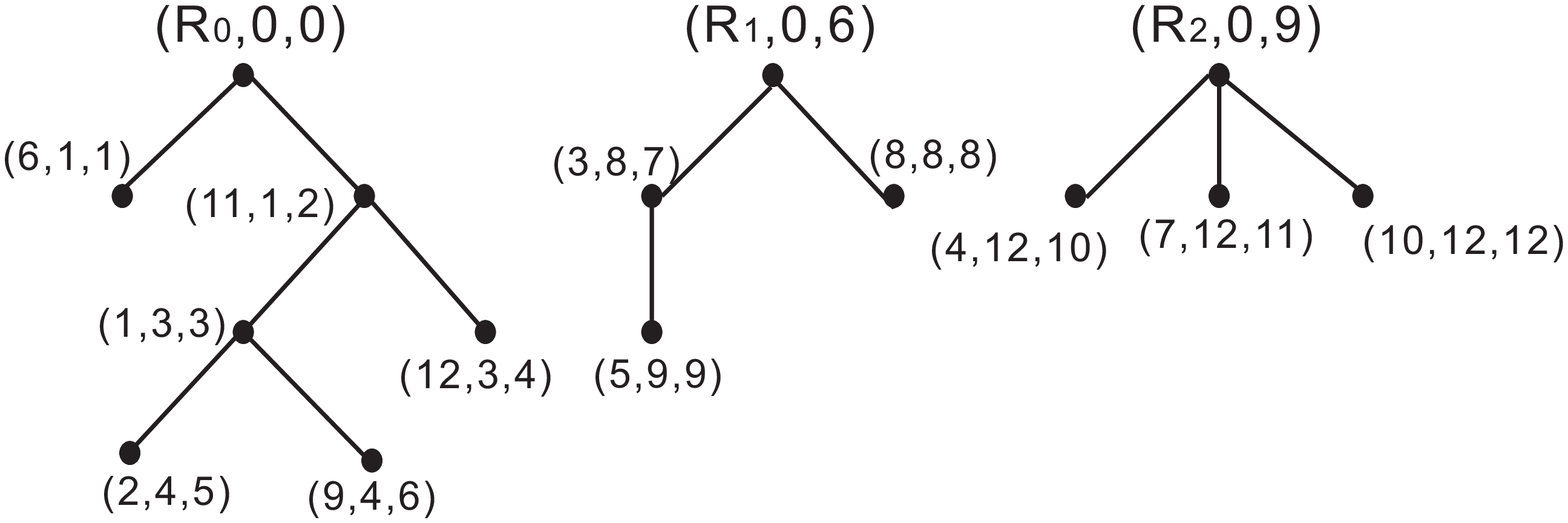}\\
Fig.\thefigure.  The corresponding forest $\phi(F'_\alpha)$
\end{center}

\addtocounter{figure}{1}
\begin{lemma}\label{bijection<m2}
Let $k$ be an integer with $k\geq 1$. For any $1\leq l\leq n+k-1$,
Let $\mathcal{A}_2=\{\alpha\in\mathcal{P}_{n,n+k; \leq n+k}^l\mid
\tau_\alpha=l-h_\alpha, n_\alpha=1,h_\alpha\leq k-1\}$ and
$\mathcal{C}_2=\{\alpha\in\mathcal{P}_{n,n+k;\leq n+k}^{l+1}\mid
l=g_\alpha\}$, then there is a bijection from $\mathcal{A}_2$ to
$\mathcal{C}_2$.
\end{lemma}
{\bf Proof.} It suffices to establish a bijection
$\psi:\mathcal{F}({\mathcal{A}_2})\rightarrow\mathcal{F}({\mathcal{C}_2})$.
 Given an $\hat{F}_\alpha\in\mathcal{F}(\mathcal{A}_2)$, let
$u=(1,l,\pi_\alpha(1))\in \hat{F}_\alpha$. Obviously,  $u\in
T_{h_\alpha}$. Let $r_{h_\alpha+1}$ be the root of the tree
$T_{h_\alpha+1}$. Delete $T_{h_\alpha}(u)$ from $T_{h_\alpha}$ and
attach $T_{h_\alpha}(u)$ to $r_{h_\alpha+1}$ so that $u$ is the
first child of $r_{h_\alpha+1}$. Then updating the second and the
third entries of all non-root vertices by Algorithm A, and the other
trees remain unchangeable, we obtain $\psi(\hat{F}_\alpha)$. The
triplet of $u$ becomes $(1,l+1,l-h_\alpha)$. Suppose
$\hat{F}_\beta=\psi(\hat{F}_\alpha)$,
 then $g_\beta=l$, hence
$\psi(\hat{F}_\alpha)\in\mathcal{F}(\mathcal{C}_2)$.

To find $\psi^{-1}$, given an
$\hat{F}_\beta\in\mathcal{F}(\mathcal{C}_2)$, let
$u=(1,l+1,\pi_\beta(1))\in \hat{F}_\beta$. Since $l=g_\beta$, we
have $h_\beta\geq 1$ and that the parent of $u$ is the root
$R_{h_\beta}$ of the tree $T_{h_\beta}$ in $\hat{F}_\beta$. Traverse
$T_{h_\beta-1}$ by breadth-first search and suppose $v$ is the last
vertex. Delete $T_{h_\beta}(u)$ from $T_{h_\beta}$ and attach
$T_{h_\beta}(u)$ to $v$. Then updating the second and the third
entries of all non-root vertices by Algorithm A, and the other trees
remain unchangeable, we obtain $\psi^{-1}(\hat{F}_\beta)$. The
triplet of $u$ becomes $(1,l,l-h_\beta+1)$. Hence,
$\psi^{-1}(\hat{F}_\beta)\in\mathcal{F}(\mathcal{A}_2)$.
\hfill$\blacksquare$\\

For example, take $n=10$, $k=2$ and $l=3$. We consider a parking
function $\alpha=(3,10,4,10,7,1,4,1,10,7)\in\mathcal{P}_{10,12;\leq
12}$. Then $h_\alpha=0$, $\tau_\alpha=3$ and $t_\alpha=1$. Observe
that $\tau_\alpha=l-h_\alpha$. On Fig.$6$ is the corresponding
forest $\psi(\hat{F}_\alpha)$, which is obtained from
$\hat{F}_\alpha-T_1(u)$ with $T_1(u)$ attached to $(R_1,0,5)$ and
with the second and third entries of the triplet updated. Sorting
the triplets of non-root vertices by the first entries, we retrieve
the corresponding parking function $(4,10,5,10,4,1,5,1,10,4)$ with
leading term $4$ from their second entries.
\begin{center}
\includegraphics[width=3.75in]{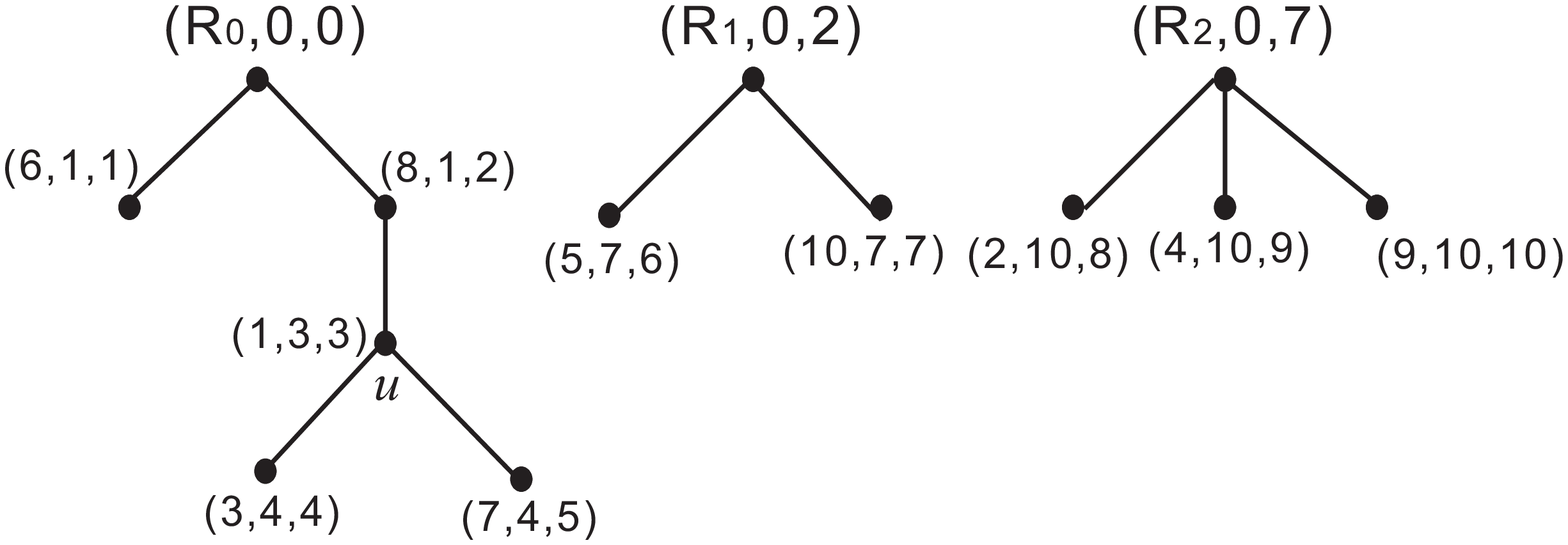}\\
Fig.\thefigure. The forest $F'_\alpha$associated with the parking
function $\alpha=(3,10,4,10,7,1,4,1,10,7)$
\end{center}
\addtocounter{figure}{1}
\begin{center}
\includegraphics[width=3.75in]{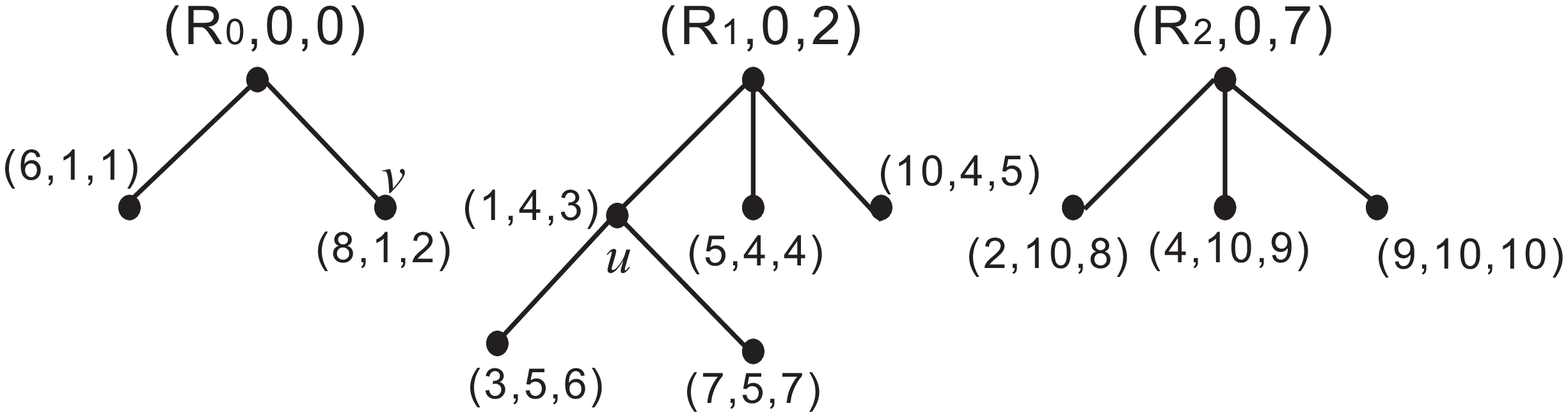}\\
Fig.\thefigure.  The corresponding forest $\phi(F'_\alpha)$
\end{center}
\addtocounter{figure}{1}
\begin{lemma}\label{bijection<m9}
Let $k\geq 0$ . For any $k+1\leq l\leq n+k$, Let
$\mathcal{A}_3=\{\alpha\in\mathcal{P}_{n,n+k;\leq n+k}^l\mid
\tau_\alpha=l-h_\alpha, n_\alpha=1, h_\alpha= k\}$ and
$\mathcal{C}_3=\{\beta\in\mathcal{P}_{n,n+k+1; \leq
n+k+1}^{n+k+1}\mid m_{\beta}=l\}$, then there is a bijection from
$\mathcal{A}_3$ to $\mathcal{C}_3$.
\end{lemma}
{\bf Proof.} For any $\alpha=(l,a_2,\ldots,a_n )\in \mathcal{A}_3$,
let $\beta=(n+k+1,a_2,\ldots,a_n)$. Let $S=\{j\mid l+1\leq a_i\leq
n+k\}$ and $T=\{j\mid m_{\alpha}+1\leq a_i\leq l-1\}$. Furthermore,
let $\alpha_S$ and $\alpha_T$ be two subsequences of $\alpha$
determined by the subscripts in $S$ and $T$, respectively. Since
$\alpha\in \mathcal{A}_3$, $\alpha_S$ and $\alpha_T$ correspond with
a parking function of length $n+k-l$ and $l-m_{\alpha}-1$,
respectively. So, $h_\alpha=k$ implies $\beta\in\mathcal{C}_3$.

Conversely, for any $\beta=(n+k+1,a_2,\ldots,a_n )\in
\mathcal{C}_3$, let $\alpha=(l,a_2,\ldots,a_n)$. $m_\beta=l$ implies
$\alpha\in \mathcal{A}_3$. This complete the proof.
\hfill$\blacksquare$
\begin{lemma}\label{lemmaA3}
Let $k\geq 0$. For any $k+1\leq l\leq n+k$, Let
$\mathcal{A}_3=\{\alpha\in\mathcal{P}_{n,n+k;\leq n+k}^l\mid
\tau_\alpha=l-h_\alpha, n_\alpha=1, h_\alpha= k\}$, then
$|\mathcal{A}_3|={n-1\choose{l-k-1}}p_{l-k-1,l-1,\leq
l-1}p_{n+k-l}$.
\end{lemma}
{\bf Proof.} Lemma \ref{bijection<m9} implies that
$|\mathcal{A}_3|=|\mathcal{C}_3|$. For any
$\alpha=(a_1,\ldots,a_n)\in\mathcal{C}_3$, let $S=\{j\mid a_i\leq
l-1\}$ and $T=[n]\setminus S$, then $|S|=l-k-1$ and $|T|=n+k+1-l$.
Furthermore, let $\alpha_S$ and $\alpha_T$ be two subsequences of
$\alpha$ determined by the subscripts in $S$ and $T$, respectively.
Obviously, $\alpha_S\in\mathcal{P}_{l-k-1,l-1,\leq l-1}$. Suppose
$\alpha_T=(b_1,\ldots,b_{n+k+1-l})$, since $l+1\leq b_j\leq n+k+1$
for any $1\leq j\leq n+k+1-l$, we have
$(b_1-l,\ldots,b_{n+k+1-l}-l)\in\mathcal{P}_{n+k+1-l}^{n+k+1-l}$.

There are ${n-1\choose{l-k-1}}$ ways to choose $l-k-1$ numbers from
$[2,n]$ for the elements in $S$. There are $p_{l-k-1,l-1,\leq l-1}$
and $p_{n+k-l}$ possibilities for $\alpha_S$ and $\alpha_T$,
respectively. Hence, we have
$|\mathcal{A}_3|=|\mathcal{C}_3|={n-1\choose{l-k-1}}p_{l-k-1,l-1,\leq
l-1}p_{n+k-l}$.\hfill$\blacksquare$

\section{Enumerations for parking functions in
$\mathcal{P}_{n,n+k;\leq n+k}^l$}

In this section, with the benefit of the triplet-labelled rooted
forests, we enumerate parking functions in the set
$\mathcal{P}_{n,n+k;\leq n+k}^l$.

\begin{theorem}\label{theoremp(n,n+k,<=n+k)1=k+1}
For any $k\geq 0$, we have $p_{n,n+k;\leq n+k}^1=p_{n,n+k;\leq
n+k}^2=\ldots =p_{n,n+k;\leq n+k}^{k+1}$.
\end{theorem}
{\bf Proof.} Let $\mathcal{A}_1$ and $\mathcal{C}_1$ be defined as
that in Lemma \ref{bijection<m1}, $\mathcal{A}_2$ and
$\mathcal{C}_2$ defined as that in Lemma \ref{bijection<m2}. Note
that $\mathcal{P}_{n,n+k;\leq
n+k}^l=\bigcup\limits_{i=1}^{2}\mathcal{A}_i$ and
$\mathcal{P}_{n,n+k;\leq
n+k}^{l+1}=\bigcup\limits_{i=1}^{2}\mathcal{C}_i$. Hence, by Lemmas
\ref{bijection<m1} and \ref{bijection<m2}, we obtain the desired
results.\hfill$\blacksquare$

\begin{theorem}\label{theoremp(n,n+k;<=n+k)^1}
(1) $p^1_{n,n+k;\leq n+k}=p_{n-1,n+k;\leq n+k}$ for any $k\geq 0$
and $n\geq 1$;

(2) $p_{n,n+k, \leq n+k}^1
=\sum\limits_{i=1}^{n}{n-1\choose{i-1}}p_{i,i+k-1,\leq
i+k-1}^1p_{n-i}$ for any $k\geq 1$.
\end{theorem}
{\bf Proof.} (1) For any
$\alpha=(1,a_2,\ldots,a_n)\in\mathcal{P}_{n,n+k;\leq n+k}^1$, let
$\psi(\alpha)=(a_2,\ldots,a_n)$, then
$\psi(\alpha)\in\mathcal{P}_{n-1,n+k;\leq n+k}$. Obviously, the
mapping $\psi$ is a bijection between the sets
$\mathcal{P}_{n,n+k;\leq n+k}^1$ and $\mathcal{P}_{n-1,n+k;\leq
n+k}$. Hence, $p^1_{n,n+k;\leq n+k}=p_{n-1,n+k;\leq n+k}$ for any
$k\geq 0$ and $n\geq 1$.

(2) For any $\alpha\in \mathcal{P}_{n,n+k,\leq n+k}^1$, we suppose
that the last parking space which don't be occupied is $i$.
Obviously, $k+1\leq i\leq n+k$. Let $S=\{j\mid a_j\leq
i-1,a_j\in\alpha\}$ and $T=[n]\setminus S$, then $|S|=i-k$,
$|T|=n+k-i$. Let $\alpha_{S}$ and $\alpha_{T}$ be two subsequences
of $\alpha$ determined by the subscripts in $S$ and $T$,
respectively. Then we have $\alpha_S\in\mathcal{P}_{i-k,i-1,\leq
i-1}^1$. Suppose $\alpha_{T}=(b_1,\ldots,b_{n+k-i})$, then
$(b_1-i,\ldots,b_{n+k-i}-i)\in \mathcal{P}_{n+k-i}$.

There are $n-1\choose{i-k-1}$ ways to choose $i-k-1$ numbers from
$[2,n]$ for elements in $S$ since $1\in S$. There are
$p_{i-k,i-1,\leq i-1}^1$ and $p_{n+k-i}$ possibilities for the
preference sets $\alpha_S$ and $\alpha_T$, respectively.

Hence, we have \begin{eqnarray*}p_{n,n+k, \leq
n+k}^1&=&\sum\limits_{i=k+1}^{n+k}{n-1\choose{i-k-1}}p_{i-k,i-1,\leq
i-1}^1p_{n+k-i}\\
&=&\sum\limits_{i=1}^{n}{n-1\choose{i-1}}p_{i,i+k-1,\leq
i+k-1}^1p_{n-i}.\end{eqnarray*}\hfill$\blacksquare$

\begin{theorem}\label{theoremrecp(n,n+k;<=n+k)^l^i+1}
Let $n\geq 1$, $k\geq 0$ and $k+1\leq l\leq n+k-1$. Then
$$p_{n,n+k;\leq n+k}^{l}-p_{n,n+k;\leq
n+k}^{l+1}={n-1\choose{l-k-1}}p_{l-k-1,l-1;\leq l-1}p_{n+k-l}.$$
\end{theorem}
{\bf Proof.} Let $\mathcal{A}=\{\alpha\in\mathcal{P}_{n,n+k;\leq
n+k}^{l}\mid h_\alpha=k, n_\alpha=1,\tau_\alpha=l-k\}$ and
$\mathcal{B}=\mathcal{P}_{n,n+k;\leq n+k}^{l}\backslash\mathcal{A}$.
By Lemmas \ref{bijection<m1} and \ref{bijection<m2}, there is a
bijection from the sets $\mathcal{B}$ to $\mathcal{P}_{n,n+k;\leq
n+k}^{l+1}$. Hence, $p_{n,n+k;\leq n+k}^{l}-p_{n,n+k;\leq
n+k}^{l+1}=|\mathcal{A}|$. Lemma \ref{lemmaA3} implies that
$$p_{n,n+k;\leq n+k}^{l}-p_{n,n+k;\leq
n+k}^{l+1}={n-1\choose{l-k-1}}p_{l-k-1,l-1;\leq l-1}p_{n+k-l}.$$
\hfill$\blacksquare$

\begin{lemma}\label{lemmap(n,n+k;<=n+k)^n+k=p(n-,n+k-1;<=n+k-1)}For any $k\geq 0$ and $n\geq 1$, we have
$p_{n,n+k;\leq n+k}^{n+k}=p_{n-1,n+k-1;\leq n+k-1}$
\end{lemma}
{\bf Proof.} For any
$\alpha=(n+k,a_2,\ldots,a_n)\in\mathcal{P}_{n,n+k;\leq n+k}^{n+k}$,
let  $\varphi(\alpha)=(a_2,\ldots,a_n)$, then
$\varphi(\alpha)\in\mathcal{P}_{n-1,n+k-1;\leq n+k-1}$. Obviously,
this is a bijection.\hfill$\blacksquare$

\section{The Bijections for $k$-flaw preference sets in
$\mathcal{P}_{n;\leq s;k}^l$} In this section, we consider
triplet-labelled rooted forests associated with preference sets in
the set $\mathcal{P}_{n;\leq s;k}^l$.

For any $1\leq s\leq n$, let $\mathcal{B}_{n,s,k}$ denote a set of
the forests $F=(T_0,\ldots,T_{k})\in\mathcal{B}_{n,k}$ such that (1)
the tree $T_k$ has at least $n+k-s+1$ vertices; (2) all the vertices
$\sigma_F^{-1}(s-k), \sigma_F^{-1}(s+1-k), \ldots, \sigma_F^{-1}(n)$
are leaves, where $\sigma_F^{-1}$ is the linear order on $[n]$ with
respect to $<_F$.
\begin{lemma}
Let $n$ and $k$ be two nonnegative integers. Suppose that $1\leq
s\leq n$ and the mapping $\phi$ is defined as that in Lemma
\ref{bijectionphi}, then $\phi$ is a bijection between the sets
$\mathcal{P}_{n;\leq s;k}$ and $\mathcal{B}_{n,s,k}$.
\end{lemma}
{\bf Proof.} For any $\alpha=(a_1,\ldots,a_n)\in\mathcal{P}_{n;\leq
s;k}$, if we add $k$ parking spaces, then $\alpha$ can be viewed as
an element in the set $\mathcal{P}_{n,n+k;\leq n+k}$. So, we may
suppose that the empty parking spaces are $m_1,\ldots, m_k$ and
$\hat{F}_\alpha=(T_0,\ldots,T_k)$ is the triple-let labeled rooted
forest associated with $\alpha$. Lemma \ref{bijectionphi} implies
that $\phi(\alpha)\in\mathcal{B}_{n,k}$. Since $s\leq n$, we have
$m_k\leq s-1$, hence, the tree $T_k$ has at least $n+k-s+1$
vertices. On the other hand, let ${\bf r}_{\alpha}$ and
$\pi_{\alpha}$ be the specification and the permutation on $[n]$ of
$\alpha$, respectively. Then the vertex $\pi^{-1}_{\alpha}(i)$ is a
leaf in $T_k$ since $r_{i+k+1}=0$ for any $i\in[s-k,n]$. So,
$\phi(\alpha)\in\mathcal{B}_{n,s,k}$.

Conversely, for any $F\in\mathcal{B}_{n,s,k}$, Lemma
\ref{bijectionphi} tells us that
$\phi^{-1}(F)=(a_1,\ldots,a_n)\in\mathcal{P}_{n,n+k;\leq n+k}$.
Suppose that all the empty parking spaces are $m_1,\ldots, m_k$,
then $m_k\leq s-1$ since $T_k$ has at least $n+k-s+1$ vertices. Let
$\sigma_F^{-1}$ be the linear order on $[n]$ with respect to $<_F$.
Furthermore, let $\hat{F}$ be the triple-let labeled rooted forest
 from $F$. Then $a_i\leq s$ for any $i\in[n]$ since the vertices
$\sigma_F^{-1}(s-k), \sigma_F^{-1}(s-k+1), \ldots, \sigma_F^{-1}(n)$
are leaves.  If we erase exact $k$ parking spaces $n+1,\ldots,n+k$,
then there are $k$ cars which can't be parked, hence,
$\phi^{-1}(F)\in\mathcal{P}_{n;\leq s;k}$. \hfill$\blacksquare$\\

The proofs of the following four lemmas is similar to Lemmas
\ref{bijection<m1}, \ref{bijection<m2}, \ref{bijection<m9} and
\ref{lemmaA3}. We just state them as follows.

\begin{lemma}\label{bijectionflaw<m1}
Suppose that $k\geq 1$, $l\geq 1$ and $l+1\leq s\leq n$. Let
$\mathcal{A}'_1=\{\alpha\in\mathcal{P}_{n; \leq s;k}^l\mid
h_\alpha\leq k-1\text{ and } \tau_\alpha>l-h_\alpha \text{ or }
n_\alpha\geq 2\}$ and $\mathcal{C}'_1=\{\alpha\in\mathcal{P}_{n;
\leq s;k}^{l+1}\mid h_\alpha\leq k-1\text{ and }l\geq g_\alpha+1\}$.
Then there is a bijection from $\mathcal{A}'_1$ to $\mathcal{C}'_1$.
\end{lemma}

\begin{lemma}\label{bijectionflaw<m2}
Suppose that  $k\geq 2$ and $1\leq l\leq s-3$. Let
$\mathcal{A}'_2=\{\alpha\in\mathcal{P}_{n; \leq s;k}^l\mid
\tau_\alpha=l-h_\alpha, n_\alpha=1,h_\alpha\leq k-2\}$ and
$\mathcal{C}'_2=\{\alpha\in\mathcal{P}_{n;\leq s;k}^{l+1}\mid
h_\alpha\leq k-1\text{ and } l=g_\alpha\}$, then there is a
bijection from $\mathcal{A}'_2$ to $\mathcal{C}'_2$.
\end{lemma}

\begin{lemma}\label{bijectionflaw<m9}
Suppose that  $k\geq 1$ and $k+1\leq l\leq s-2$. Let
$\mathcal{A}'_3=\{\alpha\in\mathcal{P}_{n;\leq s;k}^l\mid
\tau_\alpha=l-h_\alpha, n_\alpha=1, h_\alpha= k-1\}$ and
$\mathcal{C}'_3=\{\beta\in\mathcal{P}_{n; \leq s;k+1}^{s}\mid
m_{k}=l\}$, then there is a bijection from $\mathcal{A}'_3$ to
$\mathcal{C}'_3$.
\end{lemma}

\begin{lemma}\label{lemmaA3'}
Suppose $k\geq 1$ and $k+1\leq l\leq s-2$. Let
$\mathcal{A}'_3=\{\alpha\in\mathcal{P}_{n;\leq s;k}^l\mid
\tau_\alpha=l-h_\alpha, n_\alpha=1, h_\alpha= k-1\}$, then
$|\mathcal{A}'_3|={n-1\choose{l-k}}p_{l-k,l-1,\leq l-1}p_{n+k-l;\leq
s-l;1}^{s-l}$.
\end{lemma}

Using the above four lemmas, we may consider the problem of
enumerations of preferences sets in the set $\mathcal{P}_{n;\leq
s;k}^l$.

\section{Enumerations for the preference sets in
$\mathcal{P}_{n;\leq s;k}^l$}

First, we study the case with $1\leq l\leq k$.

\begin{theorem}
Suppose that $1\leq s\leq n$. For any $k\geq 1$, we have $p_{n;\leq
s;k}^1=p_{n;\leq s;k}^2=\ldots =p_{n;\leq s;k}^k$.
\end{theorem}
{\bf Proof.} Let $\mathcal{A}'_1$ and $\mathcal{C}'_1$ be defined as
that in Lemma \ref{bijectionflaw<m1}, $\mathcal{A}'_2$ and
$\mathcal{C}'_2$defined as that in Lemma \ref{bijectionflaw<m2}.
Note that $\mathcal{P}_{n;\leq
s;k}^l=\bigcup\limits_{i=1}^{2}\mathcal{A}'_i$ and
$\mathcal{P}_{n;\leq
s;k}^{l+1}=\bigcup\limits_{i=1}^{2}\mathcal{C}'_i$. Hence, by Lemmas
\ref{bijectionflaw<m1} and \ref{bijectionflaw<m2}, we obtain the
desired results.\hfill$\blacksquare$

\begin{theorem}\label{theoremp(n,<=s,k)1} Let $k\geq 1$ and $k+1\leq s\leq n$, then
\begin{eqnarray*} p_{n;\leq
s;k}^1&=&\sum\limits_{i=1}^{s-k-1}{n-1\choose{s-i-k-1}}p_{n+k-s+i;\leq
i}p_{s-k-i,s-i-1;\leq s-i-1}^1.
\end{eqnarray*}
\end{theorem}
{\bf Proof.} For any $\alpha=(a_1,\ldots,a_n)\in\mathcal{P}_{n;\leq
s; k}^1$, we suppose that the last empty parking space is $s-i$.
Obviously, $1\leq i\leq s-k-1$. Let $S=\{j\mid s-i+1\leq a_j\leq
s\}$ and $\alpha_S$ be a subsequence of $\alpha$ determined by the
subscripts in $S$, then $|S|=n+k-s+i$. Suppose
$\alpha_S=(b_1,\ldots,b_{n+k-s+i})$, then
$(b_1-s+i,\ldots,b_{n+k-s+i}-s+i)\in\mathcal{P}_{n+k-s+i;\leq i}$.
Let $T=[n]\setminus S$ and $\alpha_T$ be a subsequence of $\alpha$
determined by the subscripts in $T$, then $|T|=s-k-i$ and
$\alpha_T\in\mathcal{P}_{s-k-i,s-i-1;\leq s-i-1}^1$.

There are ${n-1\choose{s-k-i-1}}$ ways to choose $s-k-i-1$ numbers
from $[2,n]$ for the elements in $T$ since $1\in T$. There are
$p_{n+k-s+i}$ and $p_{s-k-i,s-i-1;\leq s-i-1}^1$ possibilities for
the parking function $\alpha_S$ and the preference set $\alpha_T$,
respectively. Hence, we have \begin{eqnarray*} p_{n;\leq
s;k}^1&=&\sum\limits_{i=1}^{s-k-1}{n-1\choose{s-i-k-1}}p_{n+k-s+i;\leq
i}p_{s-k-i,s-i-1;\leq s-i-1}^1.
\end{eqnarray*}\hfill$\blacksquare$

Now, we consider the case with $l=s$.

\begin{theorem}\label{thp(n,<=s,k+1)s=p(n,<=s,k)1}
 $p_{n;\leq
s;k+1}^s=p_{n;\leq s;k}^1$ for any $k\geq 1$ and $k+2\leq s\leq n$.
\end{theorem}
{\bf Proof.} It suffices to establish a bijection
$\varphi:\mathcal{F}({\mathcal{P}_{n;\leq
s;k}^1})\rightarrow\mathcal{F}({\mathcal{P}_{n;\leq s;k+1}^s})$.
 Given an $\hat{F}_\alpha=(T_0,\ldots,T_k)\in\mathcal{F}({\mathcal{P}_{n,\leq s,k}^1})$, let
$u=(1,1,\pi_\alpha(1))\in \hat{F}_\alpha$. Clearly, $u\in T_0$.
Deleting $T_{0}(u)$ from $T_{0}$, we denote $T'_0=T_0-T_0(u)$.
Suppose the number of the non-root vertices in $T'_0$ is $a$, let
the triplet of $u$ becomes $(R_1,0,a)$ and $T'_1=T_0(u)$. Let the
triplet $(R_2,0,\mu_2)$ of the root in $T_2$ become
$(R_{3},0,\mu_2-1)$ and let $T'_{3}=T_2$. For any $j\geq 2$, let the
triplet $(R_i,0,\mu_i)$ of the root in $T_i$ become
$(R_{i+1},0,\mu_i)$ and let $T'_{i+1}=T_i$. Traverse $T'_{k+1}$ by
breadth-first search and suppose $v$ is the $(s-1-k-\mu_k)$-th
vertex, attach a new vertex $w$, where the first entry of $w$ is
$1$,  as the first child $v$. Then updating the second and the third
entries of all non-root vertices by Algorithm A, we obtain
$\varphi(\hat{F}_\alpha)\in\mathcal{F}(\mathcal{P}_{n;\leq
s;k+1}^s)$.

To find $\varphi^{-1}$, given an
$\hat{F}_\beta=(T'_0,\ldots,T'_{k+1})\in\mathcal{F}(\mathcal{P}_{n;\leq
s;k+1}^s)$, let $w=(1,s,\pi_\beta(1))\in \hat{F}_\beta$. Clearly,
$w\in T'_{k+1}$. Delete $w$ from $T'_{k+1}$. Suppose that the root
of $T'_1$ is $u=(R_1,0,\mu_1)$. Let the first entry of the triplet
of $u$ become $1$ and attach $u$ to be the first child of the root
of $T'_0$. We denote the obtained tree as $T_0$. Let  the triplet
$(R_2,0,\mu_2)$ of the root in $T'_2$ become $(R_{1},0,\mu_2+1)$ and
let $T_{1}=T'_2$. For any $j\geq 3$, let the triplet $(R_i,0,\mu_i)$
of the root in $T'_i$ become $(R_{i-1},0,\mu_i)$ and let
$T'_{i-1}=T_i$. Then updating the second and the third entries of
all non-root vertices by Algorithm A, we obtain
$\varphi^{-1}(\hat{F}_\alpha)\in\mathcal{F}(\mathcal{P}_{n;\leq
s;k}^1)$.\hfill$\blacksquare$

\begin{theorem}\label{thp(n+1,<=s,k)s=p(n,<=s,k)}
Let $k\geq 0$ and $k+1\leq s\leq n$. Then $p_{n+1;\leq
s;k}^s=p_{n;\leq s;k}$.
\end{theorem}
{\bf Proof.} For any
$\beta=(s,b_1,\ldots,b_n)\in\mathcal{P}_{n+1,\leq s,k}^s$, we
consider the mapping $\varphi(\beta)=(b_1,\ldots,b_n)$. It is easy
to check that $\varphi$ if a bijection from the sets $p_{n+1;\leq
s;k}^s$ to $\mathcal{P}_{n;\leq s;k}$. Hence, $p_{n+1;\leq
s;k}^s=p_{n;\leq s;k}$. \hfill$\blacksquare$

\begin{corollary}
Let $k$, $n$ and $s$ be three integers. Suppose $k\geq 1$ and
$k+2\leq s\leq n$. Then $p_{n+1;\leq s;k}^1=p_{n;\leq s;k+1}$.
\end{corollary}
{\bf Proof.} By Theorems \ref{thp(n,<=s,k+1)s=p(n,<=s,k)1} and
\ref{thp(n+1,<=s,k)s=p(n,<=s,k)}, we immediately obtain the desired
results.\hfill$\blacksquare$

\begin{theorem}\label{theoremp(n,<=n;1)}
$p_{n;1}^n=p_{n}^2$ for any $n\geq 2$.
\end{theorem}
{\bf Proof.} It suffices to establish a bijection
$\varphi:\mathcal{F}({\mathcal{P}_n^2})\rightarrow\mathcal{F}({\mathcal{P}_{n,1}^n})$.
 Given an $\hat{T}_\alpha\in\mathcal{F}({\mathcal{P}_n^2})$, let
$u=(1,2,\pi_\alpha(1))\in \hat{T}_\alpha$ and $w=(R_0,0,0)$ be the
root of $\hat{T}_\alpha$. Suppose the number of the vertices in
$\hat{T}_\alpha(u)$ is $a$. Obviously, $a\geq 1$. Let the triplets
of $w$ and $u$ become $(R_1,0,a-1)$ and $(R_0,0,0)$, respectively.
We denote $T_0=\hat{T}_\alpha(u)$. Deleting $\hat{T}_{\alpha}(u)$
from $\hat{T}_{\alpha}$ and traversing
$\hat{T}_{\alpha}-\hat{T}_{\alpha}(u)$ by breadth-first search and
suppose $v$ is the $(n-1-a)$-th vertex, attaching a new vertex $s$,
where the first entry of $s$ is $1$,  as the first child $v$, we
denote the obtained tree as $T_1$. Then updating the second and the
third entries of all non-root vertices by Algorithm A, we obtain
$\psi(\hat{T}_\alpha)\in\mathcal{F}(\mathcal{P}_{n;1}^n)$.

To find $\varphi^{-1}$, given an
$\hat{F}_\beta\in\mathcal{F}(\mathcal{P}_{n;1}^n)$, let
$s=(1,n,\pi_\beta(1))\in \hat{F}_\beta=(T_0,T_1)$. Clearly, $s\in
T_1$. Suppose $u=(R_0,0,0)$ is the root of $T_0$, the number of
non-root vertices of $T_0$ is $a$ and the root of $T_1$ is
$w=(R_1,0,a)$. Delete the vertex $s$ from $T_{1}$. Let the first
entry of the triplets of $u$ become $1$ and the triplets of $w$
become $(R_0,0,0)$. Traverse $T_{1}$ by breadth-first search and
suppose $v$ is the first non-root vertex and attach $u$ to be the
first child of $v$. Then updating the second and the third entries
of all non-root vertices by Algorithm A, we obtain
$\varphi^{-1}(\hat{F}_\beta)$. The triplet of $u$ becomes
$(1,2,\pi(1))$. Hence,
$\varphi^{-1}(\hat{F}_\beta)\in\mathcal{F}(\mathcal{P}_n^2)$.\hfill$\blacksquare$

\begin{lemma}
For any $n\geq 2$, we have $p_{n;1}^{n}-p_{n;1}^{n-1}=p_{n-2}$
\end{lemma}
{\bf Proof.}  Let
$\mathcal{A}_1=\{\alpha\in\mathcal{P}_{n;1}^{n-1}\mid
r_{n-1}+r_n\leq 2\}$ and
$\mathcal{C}_1=\{\beta\in\mathcal{P}_{n;1}^n\mid r_{n-1}+r_n\leq
2\text{ and }m_{\beta}<n-1\}$. Let
$\mathcal{A}_2=\{\alpha\in\mathcal{P}_{n;1}^{n-1}\mid r_{n-1}=3\}$
and $\mathcal{C}_2=\{\beta\in\mathcal{P}_{n;1}^n\mid r_{n-1}=2
\text{ and }r_n=1\}$. For any
$\alpha=(n-1,a_2,\ldots,a_n)\in\mathcal{A}_1\cup\mathcal{A}_2$, we
consider the mapping
$\psi_1(\alpha)=(n,a_2,\ldots,a_n)\in\mathcal{C}_1\cup\mathcal{C}_2$.
Clearly, the mapping $\psi_1$ is a bijection from the sets
$\mathcal{A}_1\cup\mathcal{A}_2$ and
$\mathcal{C}_1\cup\mathcal{C}_2$.

Let $\mathcal{A}_3=\{\alpha\in\mathcal{P}_{n;1}^{n-1}\mid
r_n=1\text{ and }r_{n-1}=2\}$ and
$\mathcal{C}_3=\{\beta\in\mathcal{P}_{n;1}^n\mid r_{n-1}=1 \text{
and }r_n=2\}$. For any $\alpha=(a_1,\ldots,a_n)\in\mathcal{A}_3$,
obviously, $a_1=n-1$ and there are $i\neq 1$ and  $j$ such that
$a_i=n-1$ and $a_j=n$. Let $b_1=b_i=n$, $b_j=n-1$ and $b_h=a_h$ for
any $h\notin\{1,i,j\}$, then
$\beta=(b_1,\ldots,b_n)\in\mathcal{C}_3$. It is easy to obtain that
$\psi_2$ is a bijection between the sets $\mathcal{A}_3$ to
$\mathcal{C}_3$.

Let $\mathcal{A}_4=\{\alpha\in\mathcal{P}_{n;1}^{n-1}\mid
r_n=2\text{ and }r_{n-1}=1\}$ and
$\mathcal{C}_4=\{\beta\in\mathcal{P}_{n;1}^n\mid r_n=2\text{ and
}m_{\beta}=n-1\}$. Note that
$(\mathcal{A}_1,\mathcal{A}_2,\mathcal{A}_3,\mathcal{A}_4)$ and
$(\mathcal{C}_1,\mathcal{C}_2,\mathcal{C}_3,\mathcal{C}_4)$ form a
partition for $\mathcal{P}_{n;1}^{n-1}$ and $\mathcal{P}_{n;1}^{n}$,
respectively. Hence,
$p_{n;1}^{n}-p_{n;1}^{n-1}=|\mathcal{P}_{n;1}^{n}|-|\mathcal{P}_{n;1}^{n-1}|=|\mathcal{C}_4|-|\mathcal{A}_4|$.
It is easy to obtain that $|\mathcal{C}_4|=(n-1)p_{n-2}$ and
$|\mathcal{A}_4|={n-1\choose{2}}2(n-1)^{n-4}=(n-2)p_{n-2}$.
Therefore, $p_{n,1}^{n}-p_{n,1}^{n-1}=p_{n-2}$ for any $n\geq
2$.\hfill$\blacksquare$\\

To enumerate preference sets in $\mathcal{P}_{n;\leq s;k}^l$ with
$k\leq l\leq s$, we need the following lemma.

\begin{lemma}\label{hpos=k}
 Suppose $k\geq 1$, $k\leq l\leq s-1$ and
$s\leq n$. Let $\tilde{\mathcal{P}}_{n;\leq
s;k}^l=\{\alpha\in\mathcal{P}_{n;\leq s;k}^l\mid m_\alpha<l\}$. Then
$|\tilde{\mathcal{P}}_{n;\leq
s;k}^l|=\sum\limits_{i=k-1}^{l-2}{n-1\choose{i+1-k}}p_{i+1-k,i;\leq
i}p_{n+k-i-1;\leq s-i-1}^{l-i-1}.$
\end{lemma}
{\bf Proof.} For any
$\alpha=(a_1,\ldots,a_n)\in\tilde{\mathcal{P}}_{n;\leq s;k}^l$, we
suppose that the last empty parking space is $i+1$, then $k-1\leq
i\leq l-2$. Let $S=\{j\mid 1\leq a_j\leq i\}$ and $\alpha_S$ be a
subsequence of $\alpha$ determined by the subscripts in $S$, then
$|S|=i+1-k$ and $\alpha_S\in\mathcal{P}_{i+1-k,i;\leq i}$. Let
$T=[n]\setminus S$ and $\alpha_T$ be a subsequence of $\alpha$
determined by the subscripts in $T$, then $|T|=n+k-i-1$. Suppose
$\alpha_T=(b_1,\ldots,b_{n+k-i-1})$, since $i+2\leq a_j\leq s$ for
any $j\in T$, we have
$(b_1-i-1,\ldots,b_{n+k-i-1}-i-1)\in\mathcal{P}_{n+k-i-1;\leq
s-i-1}^{l-i-1}$.

There are ${n-1\choose{i+1-k}}$ ways to choose $i+1-k$ numbers from
$[2,n]$ for the elements in $S$ since $1\notin S$. There are
$p_{i+1-k,i;\leq i}$ and $p_{n+k-i-1;\leq s-i-1}^{l-i-1}$
possibilities for the preference set $\alpha_S$ and the parking
function $\alpha_T$, respectively. Hence, we have $$
|\tilde{\mathcal{P}}_{n;\leq
s;k}^l|=\sum\limits_{i=k-1}^{l-2}{n-1\choose{i+1-k}}p_{i+1-k,i;\leq
i}p_{n+k-i-1;\leq s-i-1}^{l-i-1}. $$\hfill$\blacksquare$
\begin{example}Take $n=7$,$s=6$, $l=4$ and $k=2$. By the computer search,
we have $|\tilde{\mathcal{P}}_{7;\leq 6;2}^4|=6265$. By the data in
Appendix, we find $p_{0,1;\leq 1}=1$, $p_{7;\leq 4}^2=3361$,
$p_{1,2;\leq 2}=2$ and $p_{6;\leq 3}^1=242$. Hence, it is easy to
check that $|\tilde{\mathcal{P}}_{7;\leq
6;2}^4|=\sum\limits_{i=1}^{2}{6\choose{i-1}}p_{i-1,i;\leq
i}p_{8-i;\leq 5-i}^{3-i}$.
\end{example}
\begin{theorem}\label{theoremrecurrencep(n,<=s,k)l}
Let $n,s,k,l$ be integers. Suppose $k\geq 1$, $k\leq l\leq s-1$ and
$s\leq n$. Then \begin{eqnarray*}p_{n;\leq s;k}^{l+1}-p_{n;\leq
s;k}^{l}&=&\sum\limits_{i=k-1}^{l-2}{n-1\choose{i-k+1}}p_{i-k+1,i;\leq
i}\left[ p_{n+k-i-1;\leq s-i-1}^{l-i}-p_{n+k-i-1;\leq
s-i-1}^{l-i-1}\right]\\
&&+{n-1\choose{l-k}}p_{l-k,l-1;\leq l-1}\left[ p_{n+k-l;\leq
s-l}^{1}-p_{n+k-l;\leq s-l;1}^{s-l}\right].\end{eqnarray*}
\end{theorem}
{\bf Proof.} Let $\mathcal{A}'_1$ and $\mathcal{C}'_1$ defined as
that in Lemma \ref{bijectionflaw<m1}, $\mathcal{A}'_2$ and
$\mathcal{C}'_2$ defined as that in Lemma \ref{bijectionflaw<m2},
$\mathcal{A}'_3$ defined as that in Lemma \ref{bijectionflaw<m9},
$\tilde{\mathcal{P}}_{n,\leq s,k}^{l+1}$ and
$\tilde{\mathcal{P}}_{n,\leq s,k}^{l}$ defined as that in Lemma
\ref{hpos=k}. Obviously, $\mathcal{P}_{n;\leq
s;k}^{l}=\bigcup\limits_{i=1}^{3}\mathcal{A}'_i\cup\tilde{\mathcal{P}}_{n;\leq
s;k}^l$ and  $\mathcal{P}_{n;\leq
s;k}^{l+1}=\bigcup\limits_{i=1}^{2}\mathcal{C}'_i\cup\tilde{\mathcal{P}}_{n;\leq
s;k}^{l+1}$. Hence, \begin{eqnarray*}|\mathcal{P}_{n;\leq
s;k}^{l+1}|-|\mathcal{P}_{n;\leq
s;k}^{l}|=|\tilde{\mathcal{P}}_{n;\leq
s;k}^{l+1}|-|\tilde{\mathcal{P}}_{n;\leq
s;k}^{l}|-|\mathcal{A}'_3|.\end{eqnarray*} By Lemmas \ref{lemmaA3'}
and \ref{hpos=k}, we have
\begin{eqnarray*}
&&p_{n;\leq s;k}^{l+1}-p_{n;\leq
s;k}^{l}\\
&=&\sum\limits_{i=k-1}^{l-1}{n-1\choose{i+1-k}}p_{i+1-k,i;\leq
i}p_{n+k-i-1;\leq
s-i-1}^{l-i}-\sum\limits_{i=k-1}^{l-2}{n-1\choose{i+1-k}}p_{i+1-k,i;\leq
i}p_{n+k-i-1;\leq s-i-1}^{l-i-1}\\
&&-{n-1\choose{l-k}}p_{l-k,l-1;\leq l-1}p_{n+k-l;\leq s-l;1}^{s-l}\\
&=&\sum\limits_{i=k-1}^{l-2}{n-1\choose{i-k+1}}p_{i-k+1,i;\leq
i}\left[ p_{n+k-i-1;\leq s-i-1}^{l-i}-p_{n+k-i-1;\leq
s-i-1}^{l-i-1}\right]\\
&&+{n-1\choose{l-k}}p_{l-k,l-1;\leq l-1}\left[ p_{n+k-l;\leq
s-l}^{1}-p_{n+k-l;\leq s-l;1}^{s-l}\right].
\end{eqnarray*}\hfill$\blacksquare$
\begin{example}Take $n=7$,$s=6$, $l=4$ and $k=1$. By the data in
Appendix, we find $p_{7;\leq 6;1}^5=19042$, $p_{7;\leq
6;1}^4=18602$. Hence, $p_{7;\leq 6;1}^5-p_{7;\leq 6;1}^4=440$. On
the other hand, we have $p_{0,0;\leq 0}=1$, $p_{7;\leq 5}^4=9351$,
$p_{7,\leq 5}^3=10026$, $p_{1,1,\leq 1}=1$, $p_{6,\leq 4}^3=701$,
$p_{6;\leq 4}^2= 776$, $p_{2,2;\leq 2}=3$, $p_{5;\leq 3}^2=65$,
$p_{5,\leq 3}^1=80$, $p_{3,3,\leq 3}=16$, $p_{4,\leq 2}^1=8$,
$p_{4;\leq 2;1}^2=1$. It is easy to check that $p_{7;\leq
6;1}^5-p_{7;\leq 6;1}^4={6\choose{3}}p_{3,3;\leq 3}\left[p_{4;\leq
2}^1-p_{4;\leq
2;1}^2\right]+\sum\limits_{i=0}^2{6\choose{i}}p_{i,i;\leq
i}\left[p_{7-i;\leq 5-i}^{4-i}-p_{7-i;\leq 5-i}^{3-i}\right]$.
\end{example}

\section{The generating function}
In this section, we will study some generating functions for some
sequences given in the previous sections. First, we need the
following three lemmas \cite{H2}.
\begin{lemma}{\rm \cite{H2}} Suppose that $k\geq 0$. Let $Q_k(x)=\sum\limits_{n\geq 0}\frac{p_{n,n+k;\leq
n+k}}{n!}x^n$, then
 $Q_k(x)=[P(x)]^{k+1}$.\end{lemma}

 \begin{lemma}{\rm \cite{H2}} Suppose that $k\geq 0$. Let $R_k(x)=\sum\limits_{n\geq
k}\frac{p_{n;\leq n-k}}{n!}x^n$ and $R(x,y)=\sum\limits_{k\geq
0}R_k(x)y^k$, then
$R_k(x)=P(x)\sum\limits_{i=0}^{k}\frac{(-1)^i(k+1-i)^i}{i!}x^i-\sum\limits_{i=0}^{k-1}\frac{(-1)^i(k-i)^i}{i!}x^i$
and $R(x,y)=\frac{P(x)-y}{e^{xy}-y}$.\end{lemma}

\begin{lemma}{\rm \cite{H2}} Let $F(x,y,z)=\sum\limits_{k\geq 0}\sum\limits_{s\geq k}\sum\limits_{n\geq s+1}\frac{p_{n;\leq
n-k}^{n-s}}{(n-1)!}x^ny^sz^k$, then
\begin{eqnarray*}F(x,y,z)=\frac{x}{e^{xyz}-z}\left[\frac{P(xy)(P(x)-yP(xy))}{1-y}-\frac{zP(xyz)[P(x)-yzP(xyz)]}{1-yz}\right].\end{eqnarray*}\end{lemma}

Now, for any $k\geq 1$ and $s\geq 0$, we define a generating
function $D_{k,s}(x)=\sum\limits_{n\geq s+k}\frac{p_{n;\leq
n-s;k}}{n!}x^n$.

\begin{theorem}\label{theoremgeneratingp(n;n-s;k)}Suppose $s\geq 0$ and $k\geq 1$. Let
$D_{k,s}(x)$ be the generating function for $p_{n;\leq n-s;k}$, then
$$D_{k,s}(x)=
[P(x)]^{k+1}\sum\limits_{i=0}^{k+s}\frac{(-1)^i(k+s+1-i)^i}{i!}x^i-[P(x)]^{k}\sum\limits_{i=0}^{k+s-1}\frac{(-1)^i(k+s-i)^i}{i!}x^i.$$
\end{theorem}
{\bf Proof.} By Lemma \ref{lemmap(n;<=s;k)com}, we have
\begin{eqnarray*}\sum\limits_{n\geq s+k}p_{n;\leq
n-s;k}\frac{x^n}{n!}&=&\sum\limits_{n\geq
s+k}\sum\limits_{i=1}^{n-s-k}{n\choose{s+i+k}}p_{n-s-i-k,n-s-i-1;\leq
n-s-i-1}p_{s+i+k;\leq i}\frac{x^n}{n!}\\
&=&\sum\limits_{n\geq
0}\sum\limits_{i=0}^{n}\frac{p_{n-i,n+k-i-1;\leq
n+k-i-1}}{(n-i)!}\frac{p_{s+i+k;\leq i}}{(s+i+k)!}x^{n+s+k},
\end{eqnarray*}hence,
\begin{eqnarray*}D_{k,s}(x)
&=&Q_{k-1}(x)R_{s+k}(x)\\
&=&[P(x)]^{k+1}\sum\limits_{i=0}^{k+s}\frac{(-1)^i(k+s+1-i)^i}{i!}x^i-[P(x)]^{k}\sum\limits_{i=0}^{k+s-1}\frac{(-1)^i(k+s-i)^i}{i!}x^i.
\end{eqnarray*}\hfill$\blacksquare$

\begin{corollary}
Let $D_{k}(x,y)=\sum\limits_{s\geq 0}D_{k,s}(x)y^s$ for any $k\geq
1$, then
\begin{eqnarray*}D_{k}(x,y)
=\left[\frac{P(x)}{y}\right]^k\left[R(x,y)-\sum\limits_{s=0}^{k-1}R_{s}(x)y^s\right].
\end{eqnarray*} Furthermore, let $D(x,y,z)=\sum\limits_{k\geq
1}D_{k}(x,y)z^k$, then
\begin{eqnarray*}D(x,y,z)
=\frac{zP(x)}{y-zP(x)}\left[\frac{P(x)-y}{e^{xy}-y}-\frac{(1-z)P(x)}{e^{xzP(x)}-zP(x)}\right].
\end{eqnarray*}
\end{corollary}
{\bf Proof.} By Theorem \ref{theoremgeneratingp(n;n-s;k)}, we have
\begin{eqnarray*}D_{k}(x,y)&=&\sum\limits_{s\geq 0}D_{k,s}(x)y^s\\
&=&\sum\limits_{s\geq 0}[P(x)]^kR_{s+k}(x)y^s\\
&=&\left[\frac{P(x)}{y}\right]^k\left[R(x,y)-\sum\limits_{s=0}^{k-1}R_{s}(x)y^s\right].
\end{eqnarray*}
Hence,
\begin{eqnarray*}D(x,y,z)
&=&\sum\limits_{k\geq
1}\left[\frac{P(x)}{y}\right]^k\left[R(x,y)-\sum\limits_{s=0}^{k-1}R_{s}(x)y^s\right]z^k\\
&=&\sum\limits_{k\geq
0}\left[\frac{zP(x)}{y}\right]^{k+1}R(x,y)-\sum\limits_{k\geq
0}\left[\frac{zP(x)}{y}\right]^{k+1}\sum\limits_{s=0}^{k}R_{s}(x)y^s\\
&=&\frac{zP(x)}{y-zP(x)}[R(x,y)-R(x,zP(x)]\\
&=&\frac{zP(x)}{y-zP(x)}\left[\frac{P(x)-y}{e^{xy}-y}-\frac{(1-z)P(x)}{e^{xzP(x)}-zP(x)}\right].
\end{eqnarray*}\hfill$\blacksquare$\\

For any $k\geq 1$, we define a generating function
$I_k(x)=\sum\limits_{n\geq 1}\frac{p_{n,n+k;\leq n+k}^1}{(n-1)!}x^n$
and let $I(x,y)=\sum\limits_{k\geq 0}I_{k}(x)y^k$.
\begin{lemma}Suppose that $k\geq 0$. Let $I_{k}(x)$ be the
generating function for $p_{n,n+k;\leq n+k}^1$, then $I_k(x)$
satisfies the following recurrence relation
$$I_k(x)=I_{k-1}(x)P(x)$$ for any $k\geq 1$, with the initial
condition $I_0(x)=x[P(x)]^2$. Equivalently,
$$I_{k}(x)=x[P(x)]^{k+2}.$$ Let $I(x,y)=\sum\limits_{k\geq
0}I_k(x)y^k$, then $$I(x,y)=\frac{x[P(x)^2]}{1-yP(x)}.$$
\end{lemma}
{\bf Proof.} When $k=0$, it is known that $I_0(x)=x[P(x)]^2$. When
$k\geq 1$, Theorem \ref{theoremp(n,n+k;<=n+k)^1} implies that
\begin{eqnarray*}\sum\limits_{n\geq 1}\frac{p_{n,n+k;\leq
n+k}^1}{(n-1)!}x^n=\sum\limits_{n\geq
1}\sum\limits_{i=1}^n\frac{p_{i,i+k-1;\leq
i+k-1}^1}{(i-1)!}\frac{p_{n-i}}{(n-i)!}x^n.
\end{eqnarray*}
Hence,
$$I_k(x)=I_{k-1}(x)P(x)\text{~~and~~}I_{k}(x)=x[P(x)]^{k+2}.$$ It is
easy to obtain that $$I(x,y)=\frac{x[P(x)]^2}{1-yP(x)}$$
\hfill$\blacksquare$\\

 For any $l\geq 0$ and $k\geq 0$, we define a generating
function $H_{l,k}(x)=\sum\limits_{n\geq l+1}\frac{p_{n,n+k;\leq
n+k}^{n+k-l}}{(n-1)!}x^n$
\begin{theorem}\label{theoremgeneratingp(n,n+k;<=n+k)^n+k-l} Suppose $l\geq 0$ and $k\geq 0$. Let $H_{l,k}(x)$ be
the generating function for $p_{n,n+k;\leq n+k}^{n+k-l}$, then
$H_{l,k}(x)$ satisfies the following recurrence relation
\begin{eqnarray*}H_{l,k}(x)=H_{l-1,k}(x)-\frac{p_{l,l+k;\leq
l+k}^{1}}{(l-1)!}x^l+\frac{p_{l}}{l!}x^{l+1}[P(x)]^{k+1}\end{eqnarray*}
with the initial conditions $H_{0,k}(x) =x[P(x)]^{k+1}. $
\end{theorem}
{\bf Proof.} When $l=0$, Lemma
\ref{lemmap(n,n+k;<=n+k)^n+k=p(n-,n+k-1;<=n+k-1)} tells us that
$$p_{n,n+k;\leq n+k}^{n+k}=p_{n-1,n+k-1;\leq n+k-1}$$ for any $k\geq
0$ and $n\geq 1$. Hence,
\begin{eqnarray*}H_{0,k}(x)&=&\sum\limits_{n\geq
1}\frac{p_{n,n+k;\leq n+k}^{n+k}}{(n-1)!}x^n\\
&=&\sum\limits_{n\geq 1}\frac{p_{n-1,n+k-1;\leq
n+k-1}}{(n-1)!}x^n\\
&=&x\sum\limits_{n\geq 0}\frac{p_{n,n+k;\leq n+k}}{n!}x^n\\
&=&xQ_{k}(x)\\
&=&x[P(x)]^{k+1}.
\end{eqnarray*}
Given $l\geq 1$, by Theorem \ref{theoremrecp(n,n+k;<=n+k)^l^i+1}, we
have \begin{eqnarray*}p_{n,n+k;\leq n+k}^{n+k-l}-p_{n,n+k;\leq
n+k}^{n+k-l+1}={n-1\choose{l}}p_{n-l-1,n+k-l-1;\leq
n+k-l-1}p_{l}.\end{eqnarray*} This  implies that
\begin{eqnarray*}H_{l,k}(x)&=&H_{l-1,k}(x)-\frac{p_{l,l+k;\leq l+k}^{k+1}}{(l-1)!}x^l+x^{l+1}\frac{p_{l}}{l!}Q_{k}(x)\\
&=&H_{l-1,k}(x)-\frac{p_{l,l+k;\leq
l+k}^{1}}{(l-1)!}x^l+\frac{p_{l}}{l!}x^{l+1}[P(x)]^{k+1}.\end{eqnarray*}\hfill$\blacksquare$

\begin{corollary}Let $H_{k}(x,y)=\sum\limits_{l\geq 0}H_{l,k}(x)y^l$
, then $H_k(x)$ satisfies the following equation
$$H_k(x,y)=yH_k(x,y)-xy[P(xy)]^{k+2}+x[P(x)]^{k+1}P(xy).$$
Equivalently,
$$H_k(x,y)=\frac{xP(xy)\left\{[P(x)]^{k+1}-y[P(xy)]^{k+1}\right\}}{1-y}.$$ Let
$H(x,y,z)=\sum\limits_{k\geq 0}H_{k}(x,y)z^k$, then
$$H(x,y,z)=\frac{xP(xy)}{1-y}\left[\frac{P(x)}{1-zP(x)}-\frac{yP(xy)}{1-zP(xy)}\right].$$
\end{corollary}
{\bf Proof.} By Theorem \ref{theoremgeneratingp(n,n+k;<=n+k)^n+k-l},
we have \begin{eqnarray*}\sum\limits_{l\geq
1}H_{l,k}(x)y^l=\sum\limits_{l\geq
1}H_{l-1,k}(x)y^l-\sum\limits_{l\geq 1}\frac{p_{l,l+k;\leq
l+k}^1}{(l-1)!}(xy)^l+x[P(x)]^{k+1}\sum\limits_{l\geq
1}\frac{p_l}{l!}(xy)^l.\end{eqnarray*} Therefore,
$H_k(x,y)=yH_k(x,y)-xy[P(xy)]^{k+2}+x[P(x)]^{k+1}P(xy).$
Equivalently, we have
\begin{eqnarray*}H_k(x,y)=\frac{xP(xy)\left\{[P(x)]^{k+1}-y[P(xy)]^{k+1}\right\}}{1-y}.\end{eqnarray*}
Furthermore,
$$H(x,y,z)=\frac{xP(xy)}{1-y}\left[\frac{P(x)}{1-zP(x)}-\frac{yP(xy)}{1-zP(xy)}\right].$$\hfill$\blacksquare$

Let $k\geq 1$ and $s\geq 0$ . Define a generating function
$M_{s,k}(x)=\sum\limits_{n\geq s+k+2}\frac{p^1_{n;\leq
n-s;k}}{(n-1)!}x^n$.
\begin{lemma}Suppose that $k\geq 1$ and $s\geq 0$. Let $M_{s,k}(x)$
be the generating function for $p^1_{n;\leq n-s;k}$, then
\begin{eqnarray*}M_{s,k}(x)=R_{s+k}(x)I_{k-1}(x).\end{eqnarray*}
Let $M_{k}(x,y)=\sum\limits_{s\geq 0}M_{s,k}(x)y^s$ and
$M(x,y,z)=\sum\limits_{k\geq 1}M_{k}(x,y)z^k$, then
\begin{eqnarray*}M_{k}(x,y)=I_{k-1}(x)y^{-k}\left[R(x,y)-\sum\limits_{s=0}^{k-1}R_{s}(x)y^s\right]\end{eqnarray*}
and
\begin{eqnarray*}M(x,y,z)
=\frac{xz[P(x)]^2}{y-zP(x)}\left[\frac{P(x)-y}{e^{xy}-y}-\frac{(1-z)P(x)}{e^{xzP(x)}-zP(x)}\right].\end{eqnarray*}
\end{lemma}
{\bf Proof.} For any $k\geq 1$ and $s\geq 0$, Theorem
\ref{theoremp(n,<=s,k)1} implies that
\begin{eqnarray*}p^1_{n;\leq
n-s;k}=\sum\limits_{i=1}^{n-s-k-1}{n-1\choose{s+k+i}}p_{s+k+i;\leq
i}p^{1}_{n-s-k-i,n-s-i-1;\leq n-s-i-1},
\end{eqnarray*}hence,
\begin{eqnarray*}M_{s,k}(x)&=&R_{s+k}(x)I_{k-1}(x).\end{eqnarray*}
So,
\begin{eqnarray*}M_{k}(x,y)&=&\sum\limits_{s\geq 0}R_{s+k}(x)I_{k-1}(x)y^s\\
&=&I_{k-1}(x)y^{-k}\left[R(x,y)-\sum\limits_{s=0}^{k-1}R_{s}(x)y^s\right]\end{eqnarray*}
Furthermore, we have
\begin{eqnarray*}M(x,y,z)&=&\sum\limits_{k\geq 1}I_{k-1}(x)y^{-k}[R(x,y)-\sum\limits_{s=0}^{k-1}R_{s}(x)y^s]z^k\\
&=&\frac{xz[P(x)]^2}{y-P(x)z}\left[R(x,y)-R(x,zP(x))\right]\\
&=&\frac{xz[P(x)]^2}{y-zP(x)}\left[\frac{P(x)-y}{e^{xy}-y}-\frac{(1-z)P(x)}{e^{xzP(x)}-zP(x)}\right].\end{eqnarray*}
\hfill$\blacksquare$\\

Define a generating function
$$W(x,y,z,v)=\sum\limits_{k\geq 1}\sum\limits_{s\geq
0}\sum\limits_{l\geq s}\sum\limits_{n\geq k+l}\frac{p_{n;\leq
n-s;k}^{n-l}}{(n-1)!}x^ny^lz^sv^k.$$

\begin{theorem} \begin{eqnarray*}
W(x,y,z,v)&=&\frac{xyvP(xy)}{y-1}\left\{\left[\frac{yP(xy)}{yz-vP(xy)}-\frac{P(x)}{yz-vP(x)}\right]R(xy,z)\right.\\
&&+\left.\frac{P(x)}{yz-vP(x)}R(xy,\frac{v}{y}P(x))-\frac{yP(xy)}{yz-vP(xy)}R(xy,\frac{v}{y}P(xy))\right\}\\
&&+\frac{vP(x)}{yz-vP(x)}\left[F(x,y,z)-F(x,y,\frac{v}{y}P(x))\right]\end{eqnarray*}
\end{theorem}
{\bf Proof.} First, let $W_{k,s,l}(x)=\sum\limits_{n\geq
k+l}\frac{p_{n;\leq n-s;k}^{n-l}}{(n-1)!}x^n.$ When  $l=s=0$,
Theorems \ref{theoremp(n,<=n;1)} and
\ref{thp(n,<=s,k+1)s=p(n,<=s,k)1} imply that
$$p^n_{n;k}=\left\{\begin{array}{lll}
p_n^2=p_n^1-p_{n-1}&if&k=1\\
p_{n;k-1}^1&if&k\geq 2\end{array}\right..$$
This tells us that
$$W_{k,0,0}(x)=\left\{\begin{array}{lll}
xP(x)[P(x)-1]&if&k=1\\
M_{0,k-1}(x)&if&k\geq 2\end{array}\right..$$When $l=s\geq 1$, by
Theorem \ref{thp(n+1,<=s,k)s=p(n,<=s,k)}, we have
$W_{k,s,s}(x)=xD_{k,s-1}(x)$.\\ For any $l\geq s+1$, by Theorem
\ref{theoremrecurrencep(n,<=s,k)l}, it follows that
\begin{eqnarray*}p_{n;\leq n-s;k}^{n-l+1}-p_{n;\leq n-s;k}^{n-l}&=&\sum\limits_{i=0}^{n-k-l-1}
{n-1\choose{i}}p_{i,i+k-1;\leq i+k-1}\left[p_{n-i;\leq
n-s-k-i}^{n-k-l+1-i}-p_{n-i;\leq
n-s-k-i}^{n-k-l-i}\right]\\
&&+{n-1\choose{n-k-l}}p_{n-k-l,n-l-1;\leq n-l-1}\left[p_{l+k;\leq
l-s}^1-p_{l+k-1;\leq l-s;1}\right]\\
&=&\sum\limits_{i=0}^{n-k-l} {n-1\choose{i}}p_{i,i+k-1;\leq
i+k-1}\left[p_{n-i;\leq n-s-k-i}^{n-k-l+1-i}-p_{n-i;\leq
n-s-k-i}^{n-k-l-i}\right]\\
&&-{n-1\choose{n-k-l}}p_{n-k-l,n-l-1;\leq n-l-1}p_{l+k-1;\leq
l-s;1}.\\\end{eqnarray*}
Let $F_{k,s}(x)=\sum\limits_{n\geq
s+1}\frac{p^{n-s}_{n;\leq n-k}}{(n-1)!}x^n$, then
\begin{eqnarray*}
&&W_{k,s,l-1}(x)-\frac{p^{k}_{k+l-1;\leq
k+l-1-s;k}}{(k+l-2)!}x^{k+l-1}-W_{k,s,l}(x)\\
&=&\left[-\frac{p_{l+k-1;\leq
l-s;1}}{(k+l-1)!}x^{k+l}+F_{s+k,k+l-1}(x)-F_{s+k,k+l}(x)\right]Q_{k-1}(x).
\end{eqnarray*}
Let $W_{k,s}(x,y)=\sum\limits_{l\geq s}\sum\limits_{n\geq
k+l}\frac{p_{n;\leq n-s;k}^{n-l}}{(n-1)!}x^ny^l=\sum\limits_{l\geq
s}W_{k,s,l}(x)y^l$ and $F_{k}(x,y)=\sum\limits_{s\geq
k}F_{k,s}(x)y^s$, then
\begin{eqnarray*}&&yW_{k,s}(x,y)-y^{-k+1}M_{s,k}(xy)-W_{k,s}(x,y)+W_{k,s,s}(x)y^s\\
&=&\left[-D_{1,s+k-1}(xy)y^{-k+1}x+F_{s+k}(x,y)y^{-k}(y-1)+F_{s+k,s+k}(x)y^s\right]Q_{k-1}(x).\end{eqnarray*}
Furthermore, note that $F(x,y,z)=\sum\limits_{k\geq 0}F_{k}(x,y)z^k$
and let $W_k(x,y,z)=\sum\limits_{s\geq 0}\sum\limits_{l\geq
s}\sum\limits_{n\geq k+l}\frac{p_{n;\leq
n-s;k}^{n-l}}{(n-1)!}x^ny^lz^s=\sum\limits_{s\geq
0}W_{k,s}(x,y)z^s,$ then
\begin{eqnarray*}&&(y-1)W_{k}(x,y,z)-y^{-k+1}M_{k}(xy,z)+W_{k,0,0}(x)+xyzD_{k}(x,yz)\\
&=&-\left[D_1(xy,z)-\sum\limits_{s=0}^{k-2}D_{1,s}(xy)z^s\right](yz)^{-k+1}xQ_{k-1}(x)\\
&&+\left[F(x,y,z)-\sum\limits_{s=0}^{k-1}F_{s}(x,y)z^s\right](yz)^{-k}(y-1)Q_{k-1}(x)\\
&&+x(yz)^{-k+1}\left[R(x,yz)-\sum\limits_{s=0}^{k-2}R_{s}(x)(yz)^s\right]Q_{k-1}(x)\\
&&-x^{k}\sum\limits_{s\geq k-1}\frac{p_{s;\leq
0}}{s!}(xyz)^{s-k+1}Q_{k-1}(x).\end{eqnarray*}
Since
$W(x,y,z,v)=\sum\limits_{k\geq 1}\sum\limits_{s\geq
0}\sum\limits_{l\geq s}\sum\limits_{n\geq k+l}\frac{p_{n;\leq
n-s;k}^{n-l}}{(n-1)!}x^ny^lz^sv^k=\sum\limits_{k\geq
1}W_k(x,y,z)v^k$, we have
\begin{eqnarray*}&&(y-1)W(x,y,z,v)-yM(xy,z,\frac{v}{y})+vxP(x)[R(x,vP(x))-1]+xyzD(x,yz,v)\\
&=&-\frac{xyzvP(x)}{yz-vP(x)}D_1(xy,z)+\frac{xv^2[P(x)]^2}{yz-vP(x)}D_1(xy,\frac{v}{y}P(x))\\
&&+\frac{(y-1)vP(x)}{yz-vP(x)}\left[F(x,y,z)-F(x,y,\frac{v}{y}P(x))\right]\\
&&+\frac{xyzvP(x)}{yz-vP(x)}R(x,yz)-\frac{xv^2[P(x)]^2}{yz-vP(x)}R(x,vP(x))-xvP(x)\end{eqnarray*}
Therefore,
\begin{eqnarray*}
W(x,y,z,v)&=&\frac{xyvP(xy)}{y-1}\left[\frac{yP(xy)}{yz-vP(xy)}-\frac{P(x)}{yz-vP(x)}\right]R(xy,z)\\
&&+\frac{xyvP(xy)}{y-1}\left[\frac{P(x)}{yz-vP(x)}R(xy,\frac{v}{y}P(x))-\frac{yP(xy)}{yz-vP(xy)}R(xy,\frac{v}{y}P(xy))\right]\\
&&+\frac{vP(x)}{yz-vP(x)}\left[F(x,y,z)-F(x,y,\frac{v}{y}P(x))\right]\end{eqnarray*}
\hfill$\blacksquare$

\section{Appendix}
For convenience to check the equations given in the previous
sections, by the computer search, for $n\leq 7$, we obtain the
number of $k$-flaw preference sets $\alpha=(a_1,\ldots,a_n)$ of
length $n$ satisfying $a_1=l$ and $a_i\leq s$ for all $i\in[n]$ and
list them in the following tables. Note that $p_{n;\leq s; k}^l=0$
if $l>s$ or $n\leq k$. In Table $5$, we give the values of
$p^l_{n,n+k;\leq n+k}$ for any $n\leq 5$ and $k\leq 3$.

\newpage{\tiny
$$
\begin{array}{|r|r|l|l|l|l|l|l|l|l|}
\hline
 k=0&l=1&2&3&4&5&6&7&8&p_{n,\leq s}\\
 \hline
(n,s)=(1,1)&1&0&&&&&&&1\\
 \hline
(2,1)&1&0&&&&&&&1\\
 \hline
(2,2)&2&1&0&&&&&&3\\
 \hline
(3,1)&1&0&&&&&&&1\\
 \hline
(3,2)&4&3&0&&&&&&7\\
 \hline
(3,3)&8&5&3&0&&&&&16\\
 \hline
(4,1)&1&0&&&&&&&1\\
 \hline
 (4,2)&8&7&0&&&&&&15\\
 \hline
 (4,3)&26&19&16&0&&&&&61\\
 \hline
 (4,4)&50&34&25&16&0&&&&125\\
 \hline
 (5,1)&1&0&&&&&&&1\\
 \hline
 (5,2)&16&15&0&&&&&&31\\
 \hline
 (5,3)&80&65&61&0&&&&&206\\
 \hline
 (5,4)&232&171&143&125&0&&&&671\\
 \hline
 (5,5)&432&307&243&189&125&0&&&1296\\
 \hline
 (6,1)&1&0&&&&&&&1\\
 \hline
 (6,2)&32&31&&&&&&&63\\
 \hline
 (6,3)&242&211&206&0&&&&&659\\
 \hline
 (6,4)&982&776&701&671&&&&&3130\\
 \hline
 (6,5)&2642&1971&1666&1456&1296&0&&&9031\\
 \hline
 (6,6)&4802&3506&2881&2401&1921&1296&0&&16807\\
 \hline
 (7,1)&1&0&&&&&&&1\\
 \hline
 (7,2)&64&63&0&&&&&&127\\
 \hline
 (7,3)&728&665&659&0&&&&&2052\\
 \hline
 (7,4)&4020&3361&3175&3130&0&&&&13686\\
 \hline
 (7,5)&14392&11262&10026&9351&9031&0&&&54062\\
 \hline
 (7,6)&36724&27693&23667&20922&18682&16807&0&&144495\\
 \hline
 (7,7)&65536&48729&40953&35328&30208&24583&16807&0&262144\\
 \hline
\end{array}
$$\begin{center} Table.1. $p_{n;\leq s}^l$ for $1\leq n\leq 7$\end{center}

$$
\begin{array}{|r|r|l|l|l|l|l|l|l|l|}
\hline
 k=1&l=1&2&3&4&5&6&7&8&p_{n,\leq s,1}\\
  \hline
(n,s)=(2,1)&0&0&&&&&&&0\\
 \hline
(2,2)&0&1&0&&&&&&1\\
 \hline
(3,1)&0&0&&&&&&&0\\
 \hline
(3,2)&0&1&0&&&&&&1\\
 \hline
(3,3)&1&4&5&0&&&&&10\\
 \hline
(4,1)&0&0&&&&&&&0\\
 \hline
 (4,2)&0&1&0&&&&&&1\\
 \hline
 (4,3)&1&8&10&0&&&&&19\\
 \hline
 (4,4)&13&29&31&34&0&&&&107\\
 \hline
 (5,1)&0&0&&&&&&&0\\
 \hline
 (5,2)&0&1&0&&&&&&1\\
 \hline
 (5,3)&1&616&19&0&&&&&36\\
 \hline
 (5,4)&23&84&97&107&0&&&&311\\
 \hline
 (5,5)&165&290&293&291&307&0&&&1346\\
 \hline
 (6,1)&0&0&&&&&&&0\\
 \hline
 (6,2)&0&1&0&&&&&&1\\
 \hline
 (6,3)&1&32&36&0&&&&&69\\
 \hline
 (6,4)&41&247&291&311&0&&&&890\\
 \hline
 (6,5)&436&1107&1206&1266&1346&0&&&5361\\
 \hline
 (6,6)&2341&3637&3591&3461&3381&3506&0&&19917\\
 \hline
 (7,1)&0&0&&&&&&&0\\
 \hline
 (7,2)&0&1&0&&&&&&1\\
 \hline
 (7,3)&1&64&69&0&&&&&134\\
 \hline
 (7,4)&75&734&857&890&0&&&&2556\\
 \hline
 (7,5)&1151&4281&4858&5161&5361&0&&&20812\\
 \hline
 (7,6)&8402&17433&18329&18602&19042&19917&0&&101725\\
 \hline
 (7,7)&37883&54690&53435&51008&48808&47433&48729&0&341986\\
 \hline
\end{array}
$$\begin{center} Table.2. $p_{n;\leq s;1}^l$ for $2\leq n\leq 7$\end{center}
$$
\begin{array}{|r|r|l|l|l|l|l|l|l|l|}
\hline
 k=2&l=1&2&3&4&5&6&7&8&p_{n,\leq s,2}\\
  \hline
(n,s)=(3,1)&0&0&&&&&&&0\\
 \hline
(3,2)&0&0&0&&&&&&0\\
 \hline
(3,3)&0&0&1&0&&&&&1\\
 \hline
 (4,1)&0&0&&&&&&&0\\
 \hline
 (4,2)&0&0&0&&&&&&0\\
 \hline
 (4,3)&0&0&1&0&&&&&1\\
 \hline
 (4,4)&1&1&8&13&0&&&&23\\
 \hline
 (5,1)&0&0&&&&&&&0\\
 \hline
 (5,2)&0&0&0&&&&&&0\\
 \hline
 (5,3)&0&0&1&0&&&&&1\\
 \hline
 (5,4)&1&1&16&23&0&&&&41\\
 \hline
 (5,5)&27&27&88&129&165&0&&&436\\
 \hline
 (6,1)&0&0&&&&&&&0\\
 \hline
 (6,2)&0&0&0&&&&&&0\\
 \hline
 (6,3)&0&0&1&0&&&&&1\\
 \hline
 (6,4)&1&1&32&41&0&&&&75\\
 \hline
 (6,5)&46&46&252&371&436&0&&&1151\\
 \hline
 (6,6)&581&581&1252&1656&1991&2341&0&&8402\\
 \hline
 (7,1)&0&0&&&&&&&0\\
 \hline
 (7,2)&0&0&0&&&&&&0\\
 \hline
 (7,3)&0&0&1&0&&&&&1\\
 \hline
 (7,4)&1&1&64&75&0&&&&141\\
 \hline
 (7,5)&81&81&740&1049&1151&0&&&3102\\
 \hline
 (7,6)&1442&1442&4572&6385&7627&8402&0&&29870\\
 \hline
 (7,7)&12373&12373&21404&26326&29938&33563&37883&0&173860\\
 \hline
\end{array}
$$\begin{center} Table.3. $p_{n;\leq s;2}^l$ for $3\leq n\leq 7$\end{center} $$
\begin{array}{|r|r|l|l|l|l|l|l|l|l|}
\hline
 k=3&l=1&2&3&4&5&6&7&8&p_{n,\leq s,3}\\
  \hline
 (n,s)=(4,1)&0&0&&&&&&&0\\
 \hline
 (4,2)&0&0&0&&&&&&0\\
 \hline
 (4,3)&0&0&0&0&&&&&0\\
 \hline
 (4,4)&0&0&0&1&0&&&&1\\
 \hline
 (5,1)&0&0&&&&&&&0\\
 \hline
 (5,2)&0&0&0&&&&&&0\\
 \hline
 (5,3)&0&0&0&0&&&&&0\\
 \hline
 (5,4)&0&0&0&1&0&&&&1\\
 \hline
 (5,5)&1&1&1&16&27&0&&&46\\
 \hline
 (6,1)&0&0&&&&&&&0\\
 \hline
 (6,2)&0&0&0&&&&&&0\\
 \hline
 (6,3)&0&0&0&0&&&&&0\\
 \hline
 (6,4)&0&0&0&1&0&&&&1\\
 \hline
 (6,5)&1&1&1&32&46&0&&&81\\
 \hline
 (6,6)&51&51&51&257&451&581&0&&1442\\
 \hline
 (7,1)&0&0&&&&&&&0\\
 \hline
 (7,2)&0&0&0&&&&&&0\\
 \hline
 (7,3)&0&0&0&0&&&&&0\\
 \hline
 (7,4)&0&0&0&1&0&&&&1\\
 \hline
 (7,5)&1&1&1&64&81&0&&&148\\
 \hline
 (7,6)&87&87&87&746&1241&1442&0&&3690\\
 \hline
 (7,7)&1763&1763&1763&4893&7942&10573&12373&0&41070\\
 \hline
\end{array}
$$\begin{center} Table.4. $p_{n;\leq s;3}^l$ for $4\leq n\leq 7$\end{center}

$$
\begin{array}{|r|r|l|l|l|l|l|l|l|l|}
\hline
&l=1&2&3&4&5&6&7&8&p_{n,n+k;\leq n+k}\\
 \hline
(n,k)=(1,0)&1&&&&&&&&1\\
 \hline
 (1,1)&1&1&&&&&&&2\\
 \hline
 (1,2)&1&1&1&&&&&&3\\
 \hline
 (1,3)&1&1&1&1&&&&&4\\
 \hline
 (2,0)&2&1&&&&&&&3\\
 \hline
(2,1)&3&3&2&&&&&&8\\
 \hline
(2,2)&4&4&4&3&&&&&15\\
 \hline
 (2,3)&5&5&5&5&4&&&&24\\
 \hline
 (3,0)&8&5&3&&&&&&16\\
 \hline
(3,1)&15&15&12&8&&&&&50\\
 \hline
(3,2)&24&24&24&21&15&&&&108\\
 \hline
(3,3)&35&35&35&35&32&24&&&196\\
 \hline
 (4,0)&50&34&25&16&&&&&125\\
 \hline
(4,1)&108&108&92&74&50&&&&432\\
 \hline
 (4,2)&196&196&196&180&153&108&&&1029\\
 \hline
 (4,3)&320&320&320&320&304&268&196&&2048\\
 \hline
 (5,0)&432&307&243&189&125&&&&1296\\
 \hline
 (5,1)&1029&1029&904&776&632&432&&&4802\\
 \hline
 (5,2)&2048&2048&2048&1923&1731&1461&1029&&12288\\
 \hline
 (5,3)&3645&3645&3645&3645&3520&3264&2832&2048&26244\\
 \hline
\end{array}
$$\begin{center} Table.5. $p_{n,n+k;\leq n+k}^l$ for $1\leq n\leq 5$ and $k\leq
3$\end{center}}



\begin{thebibliography}{99}
\bibitem{Cam} Peter J Cameron, Daniel Johannsen, Thomas Prellberg,
Pascal Schweitzer, Couting Defective Parking Functions,
arXiv:0803.0302v1, 3 Mar, 2008

\bibitem{cori2002} R. Cori, D. Rossin, B. Salvy, Polynomial ideals
for sandpiles and their Grobner bases. {\it Theoretical Computer
Science} {\bf 276} (2002), no. 1-2, 1-15.


\bibitem{EFL}Sen-Peng Eu,Tung-Shan Fu,Chun-Ju Lai,On the enumeration of parking
functions by leading terms, {\it Adv. in Appl. Math.} 35 (2005)
392-406

\bibitem{EFY} Sen-Peng Eu, Tung-Shan Fu, Yeong-NanYeh, Refined Chung-Feller theorems for lattice
paths, {\it J. Combin. Theory, Ser. A }112 (2005) 143-162.

\bibitem{ELY} Sen-Peng Eu, Shu-Chung Liu, Yeong-NanYeh, Taylor expansions for Catalan and Motzkin
numbers, {\it Adv. in Appl. Math.} 29 (2002) 345-357.

\bibitem{F}J. Fran\c{c}on, Acyclic and parking functions, {\it J. Combin. Theory Ser.
A} 18 (1975) 27-35.

\bibitem{FR} D. Foata, J. Riordan, Mappings of acyclic and parking functions,
{\it Aequationes Math.} 10 (1974) 10-22.

\bibitem{GK}J.D. Gilbey, L.H. Kalikow, Parking functions, valet functions and
priority queues, {\it Discrete Math.} 197/198 (1999) 351-373.

\bibitem{henderson1911}A. Henderson. The twenty-seven lines upon the cubic surface,
Cambridge University Press (1911)

\bibitem{H1} Po-Yi Huang, Jun Ma, Yeong-Nan Yeh, Ordered $k$-flaw
Preference sets, submitted.

\bibitem{H2} Po-Yi Huang, Jun Ma, Chun-Chen Yeh, Some enumerations for parking functions, submitted.

\bibitem{KY}J.P.S. Kung, C.H. Yan, Gon\v{c}arove polynomials and parking functions,
{\it J. Combin. Theory Ser. A} 102 (2003) 16¨C37.

\bibitem{PS}J. Pitman, R. Stanley, A polytope related to empirical
distributions, plane trees, parking functions, and the
associahedron, {\it Discrete Comput. Geom.} 27 (4) (2002) 603-634.

\bibitem{postnikov2004} Postnikov, A. and Shapiro, B. Trees, Parking
Functions, Syzygies, and Deformatioins of Monomial Ideals. {\it
Transactions of the American Mathematical Society} {\bf 356} (2004).

\bibitem{R} J.Riordan, Ballots and trees, {\it J.Combin. Theory} 6 (1969) 408-411.

\bibitem{SMP}M.P. Sch¨¹tzenberger, On an enumeration problem, {\it J. Combin. Theory} 4 (1968) 219-221.

\bibitem{SRP}R.P. Stanley, Hyperplane arrangements, interval orders and trees,
{\it Proc. Natl. Acad. Sci.} 93 (1996) 2620-2625.

\bibitem{SRP2}R.P. Stanley, Parking functions and non-crossing partitions, in: The
Wilf Festschrift, {\it Electron. J. Combin.} 4 (1997) R20.

\bibitem{Y1} C.H. Yan, Generalized tree inversions and k-parking functions, {\it J.
Combin. Theory Ser. A} 79 (1997) 268-280.

\bibitem{Y2}C.H. Yan, On the enumeration of generalized parking functions,
{\it Congr. Numer.} 147 (2000) 201-209.

\bibitem{Y3}C.H. Yan, Generalized parking functions, tree inversions and
multicolored graphs, {\it Adv. in Appl. Math.} 27 (2001) 641-670.
\end{thebibliography}
\end{document}